\documentclass[a4paper,12pt]{article}
\pdfoutput=1
\usepackage{amssymb}
\usepackage{float}
\usepackage{times}
\usepackage{graphicx}
\usepackage{amsmath}
\usepackage{subfig}
\usepackage{bm}
\usepackage{color}
\usepackage{natbib}

\newtheorem{Thm}{Theorem}
\newtheorem{Prop}[Thm]{Proposition}

\begin{document}
\thispagestyle{empty}

\begin{center} 
{\Huge\bf Reflected Generalized Beta Inverse Weibull Distribution:\\ \vspace{0.5 cm} definition and properties}\\
\ \\ \ \\ \ \\ \ \\

{\large\bf Ibrahim Elbatal} \\ \vspace{1 cm}
{
Institute of Statistical Studies and Research, Department of Mathematical Statistics \\
Cairo University - Egypt.   \\
\textbf{i\_elbatal@staff.cu.edu.eg }  \\ [1cm] }

{\large\bf Francesca Condino*} \hspace{1 cm} \textbf{and} \hspace{1 cm} {\large\bf Filippo Domma } \\ \vspace{1 cm}
{
Department of Economics, Statistics and Finance \\
University of Calabria  - Italy.  \\
\textbf{f.domma@unical.it} \hspace{0.3 cm}  and \hspace{0.3 cm}   \textbf{francesca.condino@unical.it}\\ [1cm] }
\end{center}
\vspace{5 cm}
{\large\bf * Corresponding author: francesca.condino@unical.it }

\newpage
\setcounter{page}{1}
\topmargin 2 cm
\begin{center}
{\Huge\bf Reflected Generalized Beta Inverse Weibull Distribution: \\ \vspace{0.5 cm} definition and properties}  \\ \ \\
\ \\ \ \\

\begin{abstract}
In this paper we study a broad class of distribution functions which is defined by means of reflected generalized beta distribution. This class includes that of Beta-generated distribution as a special case. In particular, we use this class to extend the Inverse Weibull distribution in order to obtain the Reflected Generalized Beta Inverse Weibull Distribution.For this new distribution, moments, entropy, order statistics and a reliability measure are derived. The link between the Inverse Weibull and the Dagum distribution is generalized. Then the maximum likelihood estimators of the parameters are examined and the observed Fisher information matrix provided. Finally, the usefulness of the model is illustrated by means of an application to real data.\\

\end{abstract}
\ \\
\end{center}
\vspace{1  cm} {\bf Key words}: Moments, entropy, reliability, bimodal distribution, compound distribution.  \\

\baselineskip 0.7 cm

\section{Introduction}
In the literature various techniques are commonly used to extend a family of distribution function; some of these are based on simple transformations of the distribution function or of the survival function. In particular, let $X$ be a non-negative random variable (\textit{rv}), representing the lifetime of an individual or unit, with cumulative distribution function (\textit{cdf}), survival function (\textit{sf}) and probability density function (\textit{pdf}), respectively, denoted by  $G(x;\boldsymbol{\tau})$, $\bar{G}(x;\boldsymbol{\tau})=1-G(x;\boldsymbol{\tau})$ and $g(x;\boldsymbol{\tau})$, where $\boldsymbol{\tau}\in\boldsymbol{\Xi}\subset \mathbb{R}^p$ with $p\geq1$ and $p\in\mathbb{N}$. The transformations $F(x;\boldsymbol{\eta})=\left[G(x;\boldsymbol{\tau})\right]^a$ and $\bar{F}(x;\boldsymbol{\eta})=\left[\bar{G}(x;\boldsymbol{\tau})\right]^b$ define the \textit{proportional reversed hazard model} (or Lehmann type I distribution) and \textit{proportional hazard model} (or Lehamann type II distribution), respectively. Usually, the distribution $G(x;\boldsymbol{\tau})$ is denoted as a parent or baseline distribution. In this context, \citet{Marshall2007} point out that the parameters $a$ and $b$ assume the meaning of \textit{resilience} and \textit{frailty} parameters, respectively. In this context, it is important to observe that once a resilience or frailty parameter has been introduced, the reintroduction of the same kind of parameter does not expand the family of distribution function; for example, in the case of the Weibull distribution, the introduction of the frailty parameter does not expand the distribution because the Weibull is already a proportional hazard family.\\
Further extensions of the families of distribution functions may be obtained by a combination of the two methods outlined above; a remarkable example in this direction is the Beta generated methods which, in brief, define the following \textit{pdf} $ f_{B-G}(x;\boldsymbol{\eta})=\frac{1}{B(a,b)}[G(x;\boldsymbol{\tau})]^{a-1} [1-G(x;\boldsymbol{\tau}))]^{b-1} g(x;\boldsymbol{\tau})$, where $\boldsymbol{\eta}=(a,b,\boldsymbol{\tau})$ and $B(\cdot,\cdot)$ is the Beta function. Thus, in a sense, the ideas of both frailty and resilience are embodied in this formulation of the \textit{pdf} \citep [pag.238]{Marshall2007}.

The class of Beta-generated distributions has received considerable attention in recent years. In particular, after the studies of \citet{Eugene2002} and \citet{Jones2004}, many Beta-generated distributions have been proposed. Among these, we recall the Beta-Gumbel  \citep{Nadarajah2004}, the Beta-Exponential \citep{Nadarajah2006} and the Beta-Weibull \citep{Famoye2005}. The reader is referred to \citet{Barreto-Souza2010},  \citet{Paranaiba2011} and  \citet{Domma2013} for the Beta-Generalized Exponential, the Beta-Burr XII and Beta-Dagum distribution, respectively, for more recent developments. A more complete list of Beta-generated distributions and an interesting review of this method for generating families of distributions can be found in \citet{Lee2013}. \\
Moreover, the Beta generated method can be obtained by means of a transformation of a random variable, \textit{i.e.} if $U\sim Be(a,b)$ then the \textit{pdf} of $X=G^{-1}(U;\boldsymbol{\eta})$ is $ f_{B-G}(x;\boldsymbol{\eta})$. \\
Recently, in order to increase the flexibility of the generator, \citet{Alexander2012} used the Generalized Beta of first type (GB1), first proposed by \citet{McDonald1984}, instead of the classic Beta function as generator. The \textit{pdf} of a GB1 is
\begin{equation}
f(u;a,b,c)=\frac{c}{B(a,b)}u^{ac-1}(1-u^{c})^{b-1}
\label{pdfGB1}
\end{equation}
for $u\in(0,1)$ and $a,b,c>0$. With this generator, the random variable $X$ is said to have a generalized beta of first type (GB1) generated distribution, with \textit{pdf} and \textit{cdf}, respectively, given by
\begin{equation}
f_{GB1-G}(x;\boldsymbol{\eta})=\frac{c}{B(a,b)}[G(x;\boldsymbol{\tau})]^{ac-1}(1-[G(x;\boldsymbol{\tau})]^{c})^{b-1}g(x;\boldsymbol{\tau})
\label{pdfGB1-G}
\end{equation}
and
\begin{equation}
F_{GB1-G}(x;\boldsymbol{\eta})=\frac{1}{B(a,b)}\int^{[G(x;\boldsymbol{\tau})]^{c}}_{0}u^{a-1}(1-u)^{b-1}du
\label{dfGB1-G}
\end{equation}
where $\boldsymbol{\eta}=(a,b,c,\boldsymbol{\tau})$. \\
In effect, as discussed in \citet{Alexander2012}, (\ref{dfGB1-G}) increases the flexibility of the generator for the presence of parameter $c$. Nevertheless, there are cases where this does not occur; in particular, when $G(x;\boldsymbol{\tau})$ is already a reversed proportional hazard distribution. In
such cases, to give effect to the expansion of the distribution function it is possible to refer to the reflected version of (\ref{pdfGB1}) obtained by transformation $v=1-u$, with \textit{pdf} and \textit{cdf}, respectively, given by
\begin{equation}
f(v;a,b,c)=\frac{c}{B(a,b)}(1-v)^{ac-1}(1-(1-v)^{c})^{b-1}
\label{pdfrGB1}
\end{equation}
\begin{equation}
F(v;a,b,c)=1-\frac{c}{B(a,b)} \int^{1-v}_{0} w^{ac-1}(1-w^{c})^{b-1}dw=1-\frac{B((1-v)^c;a,b)}{B(a,b)}
\label{dfrGB1}
\end{equation}
for $v\in(0,1)$, $a,b,c>0$, where $B(z;a,b)$ is the incomplete beta function. Putting $V=G(X;\boldsymbol{\tau})$, the random variable $X$ is said to have a reflected generalized beta of first type (rGB1) generated distribution, with \textit{pdf} and \textit{cdf}, respectively, given by
\begin{equation}
f_{rGB1-G}(x;\boldsymbol{\eta})=\frac{c}{B(a,b)}[1-G(x;\boldsymbol{\tau})]^{ac-1}(1-[1-G(x;\boldsymbol{\tau})]^{c})^{b-1}g(x;\boldsymbol{\tau})
\label{pdfrGB1-G}
\end{equation}
and
\begin{eqnarray}
F_{rGB1-G}(x;\boldsymbol{\eta})&=&1-\frac{c}{B(a,b)}\int^{[1-G(x;\boldsymbol{\tau})]}_{0}v^{ac-1}(1-v^c)^{b-1}dv= \nonumber \\
&=& 1-\frac{1}{B(a,b)}\int^{[1-G(x;\boldsymbol{\tau})]^{c}}_{0}w^{a-1}(1-w)^{b-1}dw= \nonumber \\
&=&1-\frac{B\left([1-G(x;\boldsymbol{\tau})]^{c};a,b\right)}{B(a,b)}.
\label{dfrGB1-G}
\end{eqnarray}
The class of generated distributions defined in (\ref{pdfrGB1-G}) is called \textit{Beta-Exponential-X family} by \citet{Lee2013}.
Clearly, the baseline distribution $G(x;\boldsymbol{\tau})$ is a special case of (\ref{dfrGB1-G}) when $a=b=c=1$. If $c=1$ we obtain the Beta-generated distribution. Moreover, the exponentiated generalized class of distributions proposed by \citet{Cordeiro2013} and \citet{Cordeiro2013b} turns out to be a special case of $F_{rGB1-G}(x;\boldsymbol{\eta}) $ setting $a=1$. We thus observe that the latter class of distributions can be thought of as reflected Kumaraswamy-generated distributions. \\
For $b>0$ real non-integer, using the power series representation \\
\begin{eqnarray}
(1-[1-G(x;\boldsymbol{\tau})]^{c})^{b-1}=\sum ^{\infty}_{j_{1}=0} p_{j_{1},b} [1-G(x;\boldsymbol{\tau})]^{cj_{1}}
\end{eqnarray}
where $p_{j_{1},b}=\frac{(-1)^{j_{1}} \Gamma(b) }{\Gamma(b-j_{1}) \Gamma(j_{1}+1)}$, we can write the \textit{pdf} (\ref{pdfrGB1-G}) as
\begin{equation}
f_{rGB1-G}(x;\boldsymbol{\eta})=\frac{c}{B(a,b)} \sum ^{\infty}_{j_{1}=0} p_{j_{1},b}
[1-G(x;\boldsymbol{\tau})]^{c(a+j_{1})-1}g(x;\boldsymbol{\tau})
\label{pdfrGB1-G_1}
\end{equation}
with \textit{cdf}
\begin{equation}
F_{rGB1-G}(x;\boldsymbol{\eta})=\frac{1}{B(a,b)} \sum ^{\infty}_{j_{1}=0} \frac{p_{j_{1},b}}{(a+j_1)}
\left\{1- [1-G(x;\boldsymbol{\tau})]^{c(a+j_{1})} \right\}.
\label{dfrGB1-G_1}
\end{equation}

The expansions of the \textit{pdf} shown in (\ref{pdfrGB1-G_1}) can be used to determine certain properties of \textit{rBG1-G} distribution as, for example, the \textit{k-th} moment,
\begin{eqnarray}
E\left(X^k\right) &=& \frac{c}{B(a,b)} \sum ^{\infty}_{j_{1}=0} p_{j_{1},b} \int ^{\infty}_{0} x^k
[1-G(x;\boldsymbol{\tau})]^{c(a+j_{1})-1}g(x;\boldsymbol{\tau}) dx = \nonumber \\
&& \frac{c}{B(a,b)} \sum ^{\infty}_{j_{1}=0} p_{j_{1},b} E_{G}\left\{X^k [1-G(X;\boldsymbol{\tau})]^{c(a+j_{1})-1}\right\}
\label{MrGB1-G_1}
\end{eqnarray}
where $E_G(.)$ denotes the expectation with respect to baseline distribution $G(X;\boldsymbol{\tau})$. An alternative expression for the \textit{k-th} moment is
 \begin{eqnarray}
E\left(X^k\right) = \frac{c}{B(a,b)} \sum ^{\infty}_{j_{1}=0} p_{j_{1},b} \sum ^{\infty}_{j_{2}=0} p_{j_{2},c(a+j_{1})}
 E_{G}\left\{X^k [G(X;\boldsymbol{\tau})]^{j_{2}}\right\}.
\label{MrGB1-G_2}
\end{eqnarray}
In this paper, we use the generator defined by (\ref{dfrGB1-G}) to extend the Inverse Weibull distribution. 
The rest of the paper is organized as follows. In \textit{Section 2}, we define the Reflected Generalized Beta of first type Inverse Weibull distribution and derive the expansions for its distribution and probability density functions. The quantile, the moments, the entropy, the order statistics and a measure of reliability are derived in \textit{Section 3}. \textit{In Section 4}, the link between the Inverse Weibull distribution and the Dagum distribution is generalized. The maximum likelihood estimation is discussed in \textit{Section 5}. An application on real data set is reported in \textit{Section 6}. Finally, some of the mathematical results are reported in the Appendix.

\section{Reflected Generalized Beta of first type Inverse Weibull Distribution (\textit{rGB1-IWei}) } 

\begin{figure}[ht]
\centering
\subfloat[]
{\includegraphics[scale=0.4]
{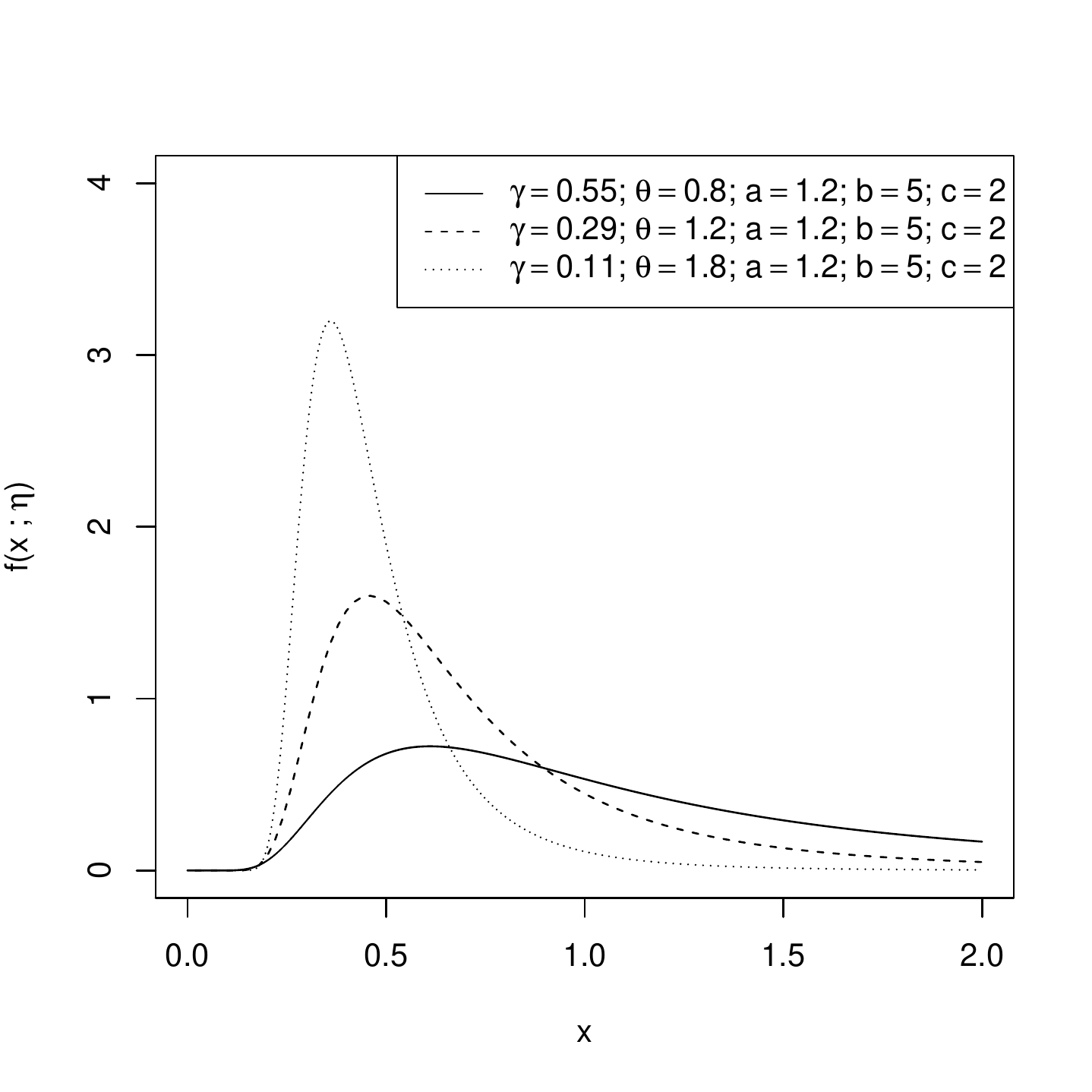} \label{dens_theta}}
\subfloat[]
{\includegraphics[scale=0.4]
{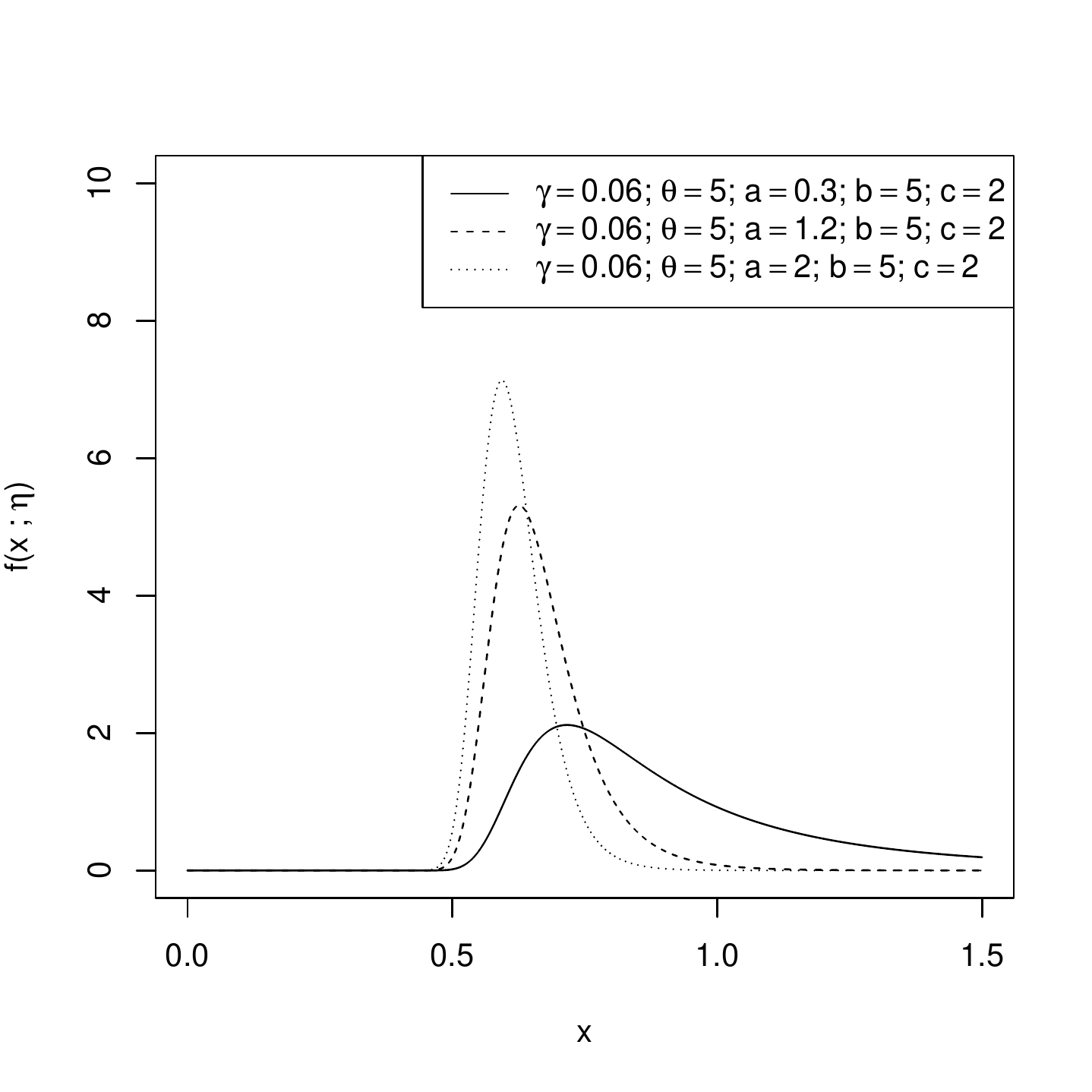}\label{dens_a}} \\
\subfloat[]
{\includegraphics[scale=0.4]
{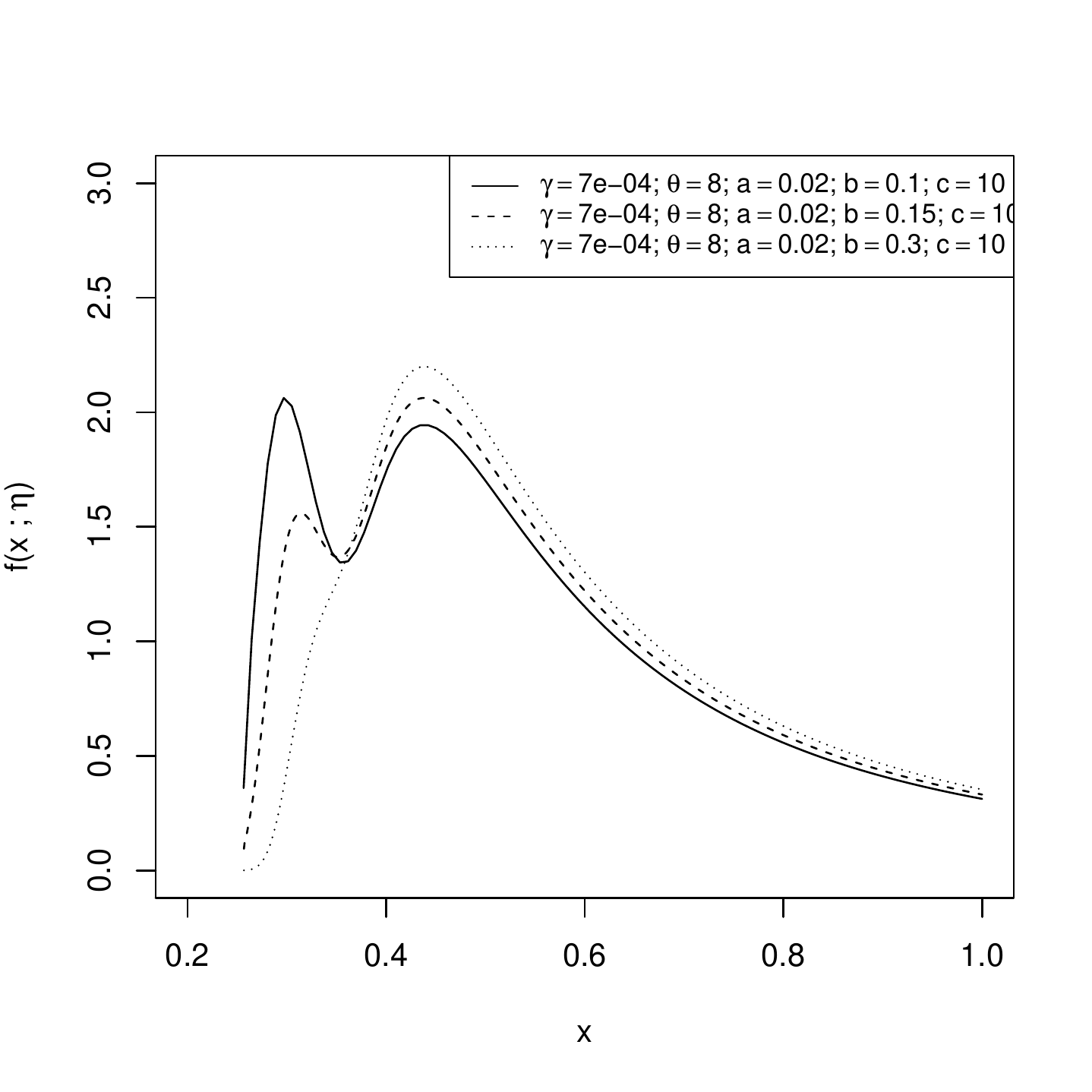}\label{dens_bimod}}
\subfloat[]
{\includegraphics[scale=0.4]
{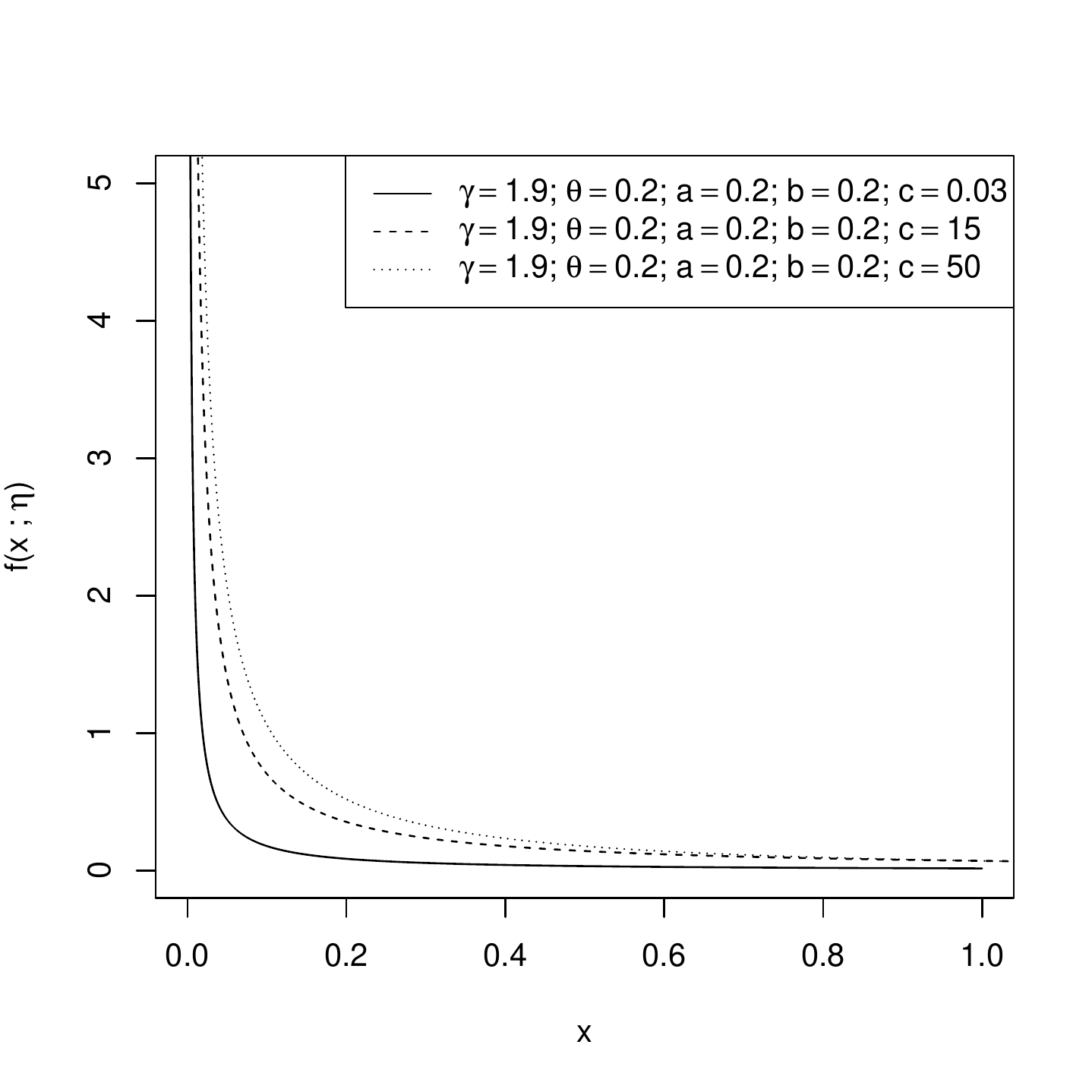}\label{dens_decr}}
\caption{\textit{rGB1-IWei} density for certain values of the parameters.}
\label{fig:dens}
\end{figure}

Several authors have highlighted the fact that there are a considerable number of contexts in which the inverse Weibull (\textit{IWei}) distribution may be an appropriate model,  mainly because the empirical hazard rate is unimodal (see, for example, \citet{Kundu2010}; \citet{Singh2013}).  Recently, \citet{Erto2013} has shown certain peculiar properties of the \textit{IWei} distribution; in particular, he formulated three real and typical degenerative mechanisms which lead exactly to the \textit{IWei} random variable. 
The cumulative distribution function and probability density function of the \textit{IWei} random variable, respectively, are:
\begin{equation}
  G_{IWei}(x;\boldsymbol{\tau})=e^{-\gamma x^{-\theta}}
\label{GIWei}
\end{equation}
and 
\begin{equation}
g_{IWei}(x;\boldsymbol{\tau})=  \theta \gamma x^{-\theta-1} e^{- \gamma x^{-\theta}},
\label{pdfGIWei}
\end{equation}
where $\boldsymbol{\tau}=(\gamma, \theta)$,  with \textit{k-th} moment about zero given by
\begin{equation}
E(X^k)= \gamma^{\frac{k}{\theta}}\Gamma \left(1-\frac{k}{\theta}\right)
\label{M_GIWei}
\end{equation}
where $\Gamma(\cdot)$ is the Gamma function. \\
Recently, \citet{Gusmao2011} have tried to generalize the \textit{IWei} distribution by introducing the \textit{resilence} parameter in (\ref{GIWei}),  \textit{i.e.} $G(x)=\left[e^{-\gamma x^{-\theta}}\right]^{\epsilon}$, 
but, as highlighted by \cite{Jones2012}, the proposed distribution is only a reparametrization of $G_{IWei}(x;\boldsymbol{\tau})$. In other words, they do not expand the class of distributions because the authors introduce a \textit{resilence} parameter in a \textit{reversed proportional hazard model}.
In order to extend the \textit{IWei} distribution, given that this distribution is a \textit{reversed hazard rate model}, we use the reflected Generalized Beta generator defined in (\ref{dfrGB1-G}). So, the \textit{rGB1-IWei} cumulative distribution function is defined by inserting (\ref{GIWei}) in (\ref{dfrGB1-G}), \textit{i.e.}
\begin{eqnarray}
F_{rGB1-IWei}(x;\boldsymbol{\eta})&=& 1-\frac{1}{B(a,b)}\int^{[1-G_{IWei}(x;\boldsymbol{\tau})]^{c}}_{0}v^{a-1}(1-v)^{b-1}dv \nonumber \\
&=& 1- I_{([1-G_{IWei}(x;\boldsymbol{\tau})]^{c})}(a,b)
\label{dfrGB1-GIWei}
\end{eqnarray}
where $\boldsymbol{\eta}=(\gamma, \theta, a, b, c)$ and $I_{G}(a,b)$ is the incomplete beta function ratio. The corresponding \textit{pdf} of the new distribution can be written as
\begin{eqnarray}
f_{rGB1-IWei}(x;\boldsymbol{\eta})= 
\frac{c \gamma \theta  x^{-\theta-1} e^{- \gamma x^{-\theta}}  }{B(a,b)}
\left[1-e^{- \gamma x^{-\theta}} \right]^{ac-1} 
 \left(1- [1-e^{- \gamma x^{-\theta}}]^{c}  \right) ^{b-1}.
\label{pdfrGB1-GIWei}
\end{eqnarray}
Fig.\ref{fig:dens} shows various behaviours of the density  $f_{rGB1-IWei}(x;\boldsymbol{\eta} )$ ,  including the bimodal one, obtained for different values of parameters $\gamma$, $\theta$, $a$, $b$ and $c$. 

\begin{figure}[ht]
\centering
\subfloat[]
{\includegraphics[scale=0.4]
{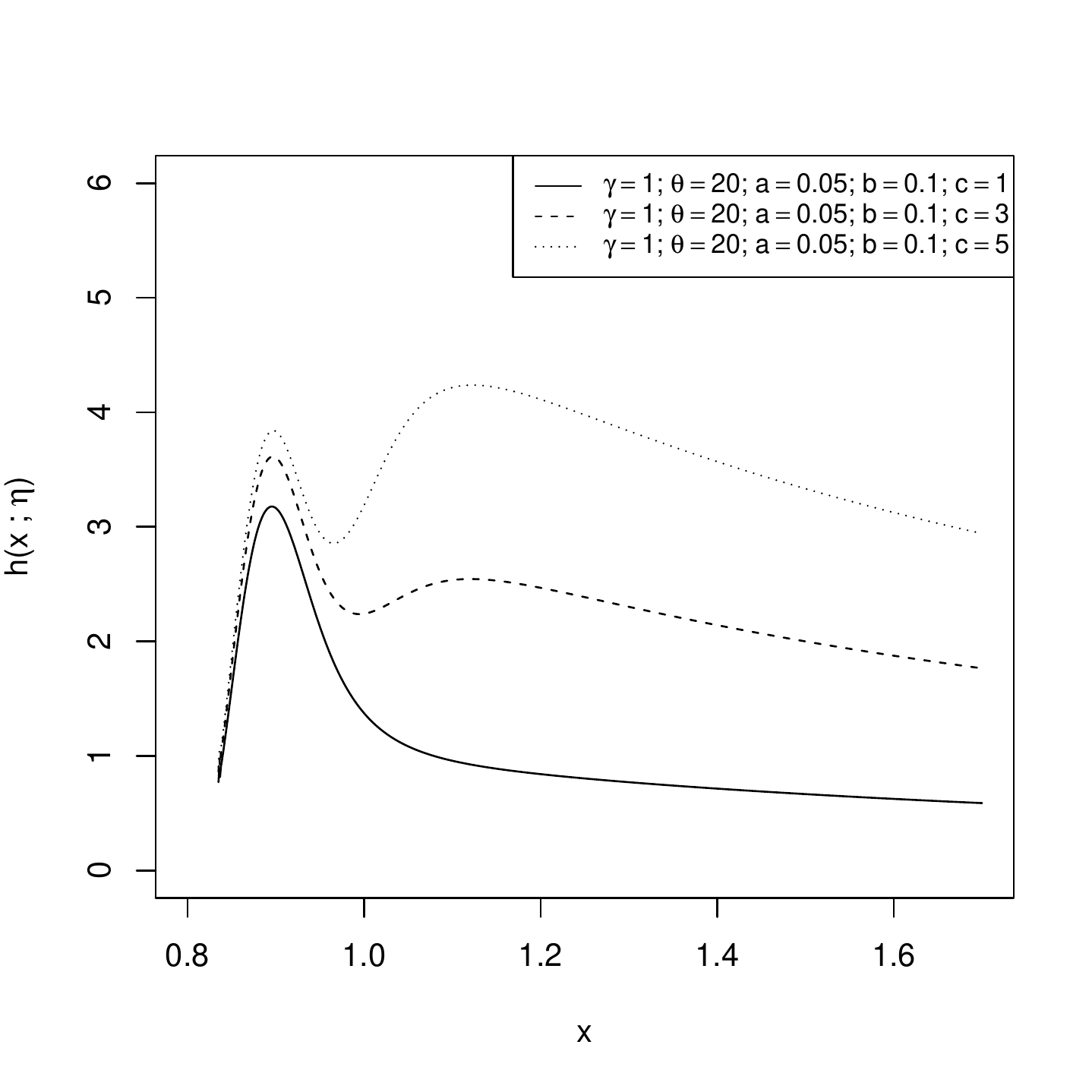} \label{hz_max&max}}
\subfloat[]
{\includegraphics[scale=0.4]
{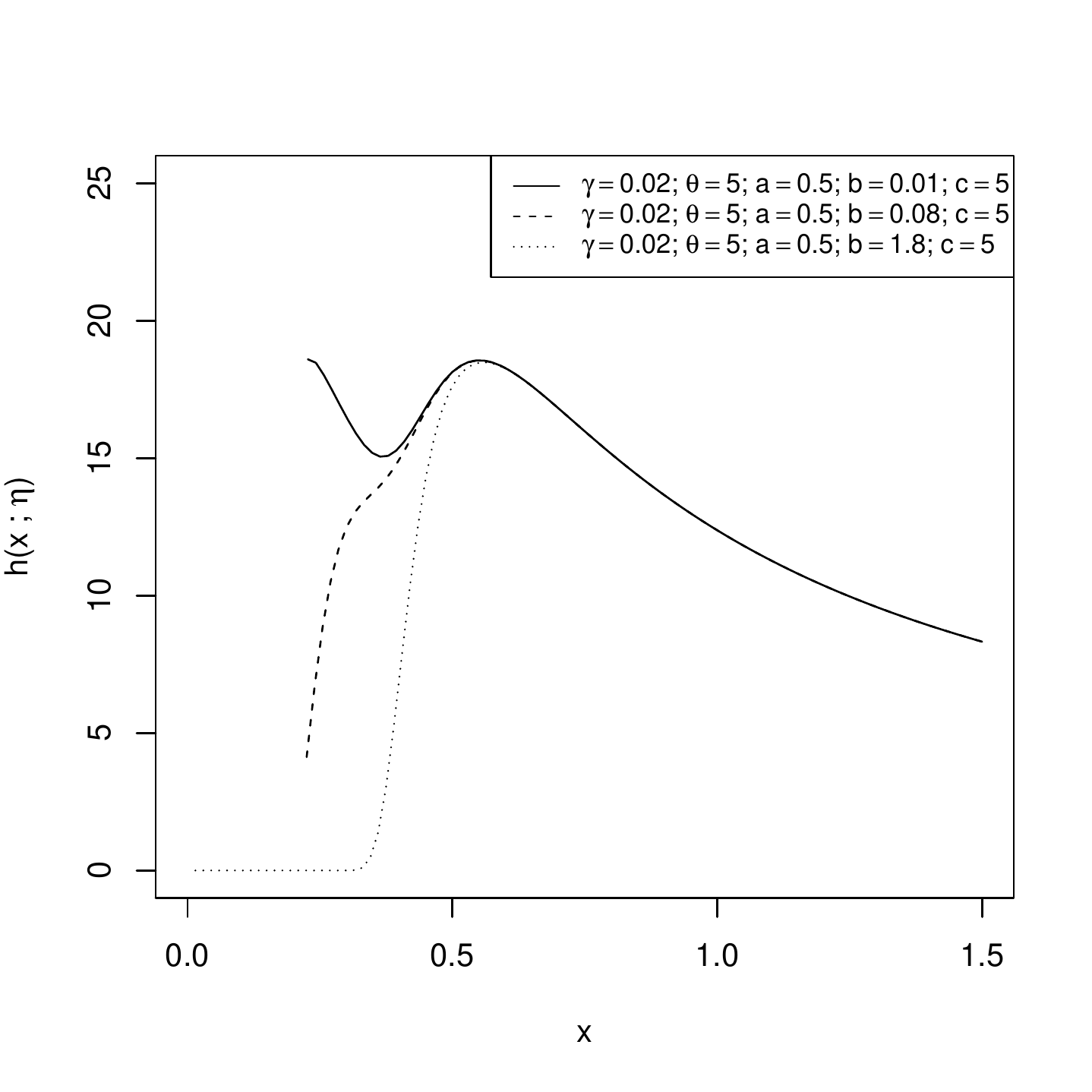}\label{hz_min&max}} \\
\subfloat[]
{\includegraphics[scale=0.4]
{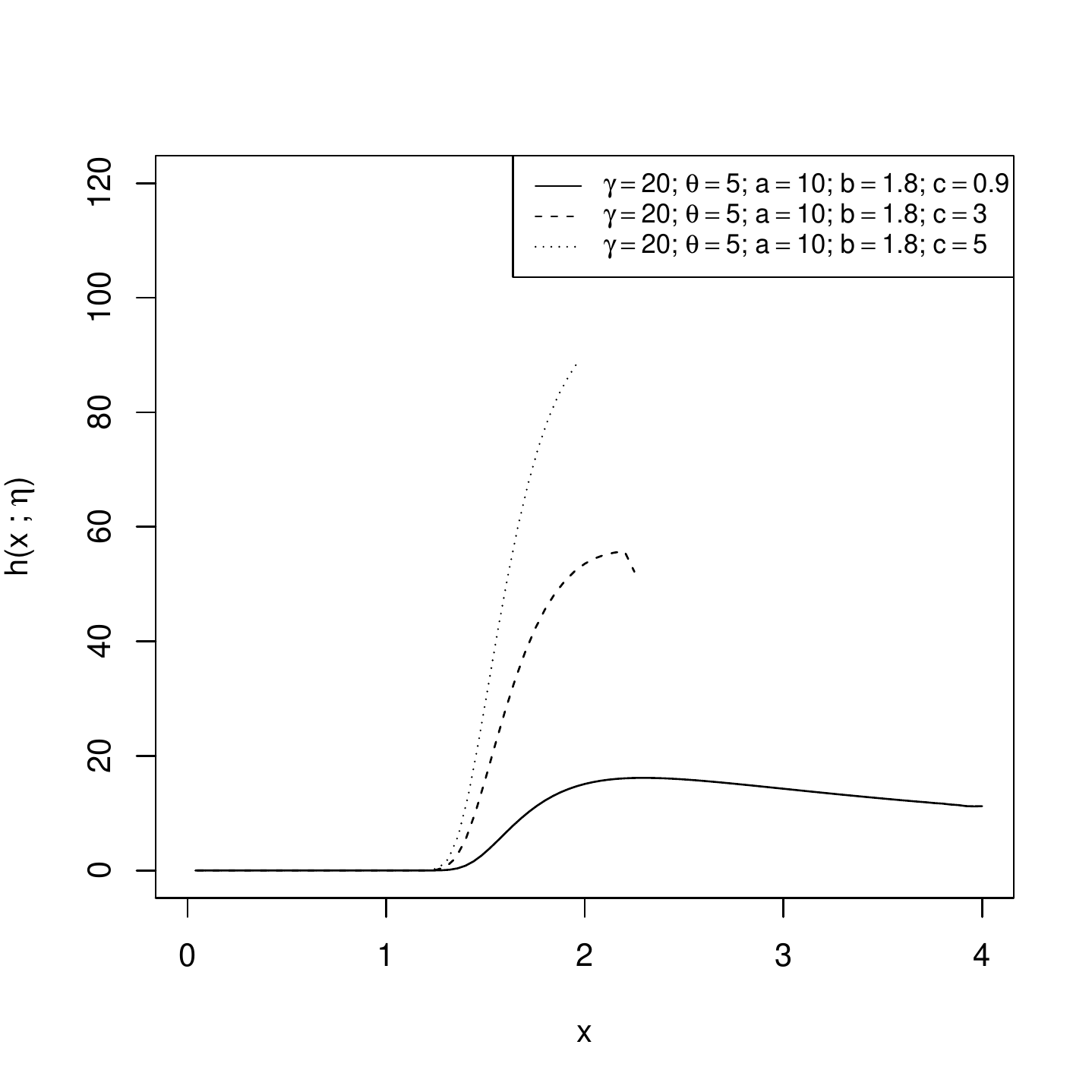}\label{hz_ifr}}
\subfloat[]
{\includegraphics[scale=0.4]
{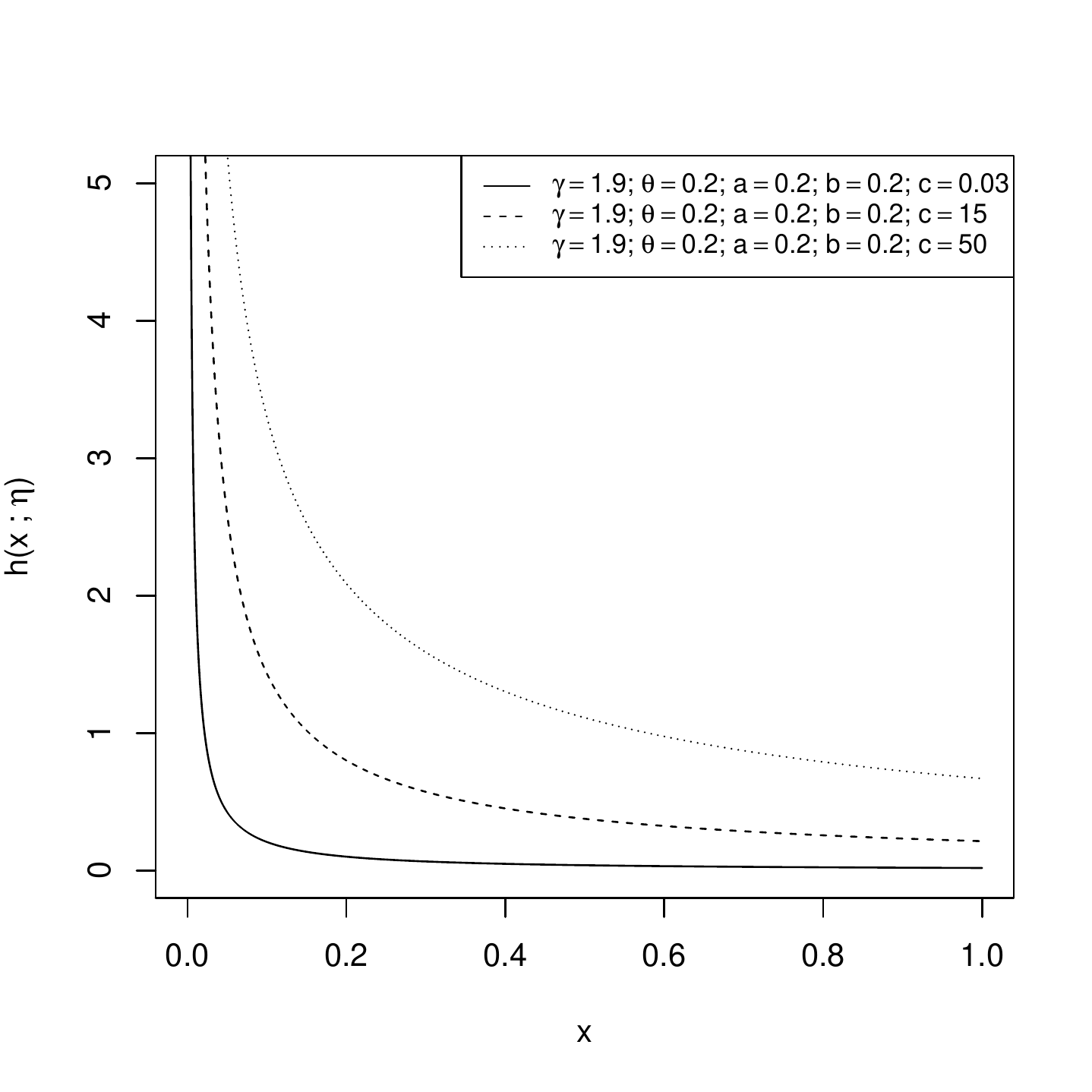}\label{hz_dfr}}
\caption{\textit{rGB1-IWei} hazard rate for certain values of the parameters.}
\label{fig:hz}
\end{figure}

The hazard rate for \textit{rGB1-IWei} distribution, given by $h_{rGB1-IWei}(x;\boldsymbol{\eta})= \frac{f_{rGB1-IWei}(x;\boldsymbol{\eta})}{I_{([1-G_{IWei}(x;\boldsymbol{\tau})]^{c})}(a,b)}$, is a complex functions of $x$. However,
from the plots reported in  Fig. \ref{fig:hz}, we can state that
$h_{rGB1-IWei}(x;\boldsymbol{\eta})$ is much more flexible than the hazard rate of the \textit{IWei} distribution. In particular, by denoting the strictly increasing (decreasing) failure
rate as IFR (DFR) and the hazard rate with a minimum or a maximum as Bathtub (BT) and Upside-down Bathtub (UBT), respectively, we can observe that the \textit{rGB1-IWei} model enables us to obtain IFR and DFR hazard rates or hazard rates with a non-monotonic behaviour, such as UBT, BT-UBT or UBT-BT-UBT hazard rate.
This wide range of different behaviours of the hazard rate makes the \textit{rGB1-IWei} distribution a suitable model for reliability theory and survival analysis.

	\subsection{Expansions for the Cumulative Distribution Function and for the Probability Density Function}
In this subsection, we present some representations of \textit{cdf} and \textsl{pdf} of \textit{rGB1-IWei} in terms of infinite sums.
Using the relations (\ref{pdfrGB1-G_1}), (\ref{GIWei}) and (\ref{pdfGIWei}), the \textit{pdf} of \textit{rGB1-IWei} can be written as
\begin{eqnarray}
f_{rGB1-IWei}(x;\boldsymbol{\eta})=\frac{c}{B(a,b)}\sum ^{\infty}_{j=0} p_{j,b}[1-G_{IWei}(x;\boldsymbol{\tau})]^{c(a+j)-1}g_{IWei}(x;\boldsymbol{\tau})
\end{eqnarray}
with \textit{cdf} given by
\begin{eqnarray}
F_{rGB1-IWei}(x;\boldsymbol{\eta})&=& \frac{1}{B(a,b)} \sum ^{\infty}_{j=0} \frac{p_{j,b}}{(a+j)} \left\{1-[1-G_{IWei}(x;\boldsymbol{\tau})]^{c(a+j)}\right\} \nonumber \\
&=& \frac{1}{B(a,b)} \sum ^{\infty}_{j=0} \frac{p_{j,b}}{(a+j)} \left\{1- [1-e^{- \gamma x^{-\theta}}]^{c(a+j)} \right\}
\end{eqnarray}

Another expansion of $f_{rGB1-IWei}(x;\boldsymbol{\eta})$ may be derived from the following relation between \textit{pdf} and \textit{cdf} of a \textit{IWei} distribution, \textit{i.e.} $g_{IWei}(x;\boldsymbol{\tau})= \gamma \theta x^{-\theta-1} G_{IWei}(x;\boldsymbol{\tau})$. After algebra, we can write $f_{rGB1-IWei}(x;\boldsymbol{\eta})$ in term of an infinite sum of \textit{IWei} density functions, \textit{i.e.}
\begin{eqnarray}
f_{rGB1-IWei}(x;\boldsymbol{\eta})=\frac{c}{B(a,b)}\sum ^{\infty}_{j_{1}=0} p_{j_{1},b}
\left(\sum ^{\infty}_{j_{2}=0}  \frac{p_{j_{2},c(a+j_{1})}}{(j_{2}+1)}      g_{IWei}(x; \gamma(j_{2}+1),\theta)\right).
\label{pdfrGB1-GIWei_2}
\end{eqnarray}
The linear combination (\ref{pdfrGB1-GIWei_2}) enables us to obtain some mathematical properties of the \textit{rGB1-IWei} distribution directly from those of the \textit{IWei} distribution such as the moments, moment generating function and reliability.

\section{Statistical Properties}
This section is devoted to studying the statistical properties of the \textit{rGB1-IWei} distribution, specifically quantile function, moments and moment generating function, entropy, order statistics and reliability.
	\subsection{Quantile of \textit{rGB1-IWei}}
The quantile $x_{q}$ of the \textit{rGB1-IWei} distribution can be easily obtained considering that%

\begin{equation*}
q=F_{rGB1-IWei}(X\leq x_q;\bm{\eta})=1-I_{([1-G_{IWei}(x_q;\bm{\tau})]^c)}(a,b)
\end{equation*}
from which
\begin{equation*}
[1-G_{IWei}(x_q;\bm{\tau})]^c=[1-e^{-\gamma  (x_q)^{\theta}}]^c=I^{-1}(1-q;a,b)
\end{equation*}
and
\begin{equation}
x_q=\left\{- \frac{1}{\gamma} \ln\left(1-[I^{-1}(1-q;a,b)]^{1/c}\right)\right\}^{-1/\theta}.
\label{xq}
\end{equation}

The median can be derived from (\ref{xq}) be setting $q=\frac{1}{2}$, that is,
\begin{equation}
Med_{(rGB1-IWei)}= \left\{- \frac{1}{\gamma} \ln\left(1-[I^{-1}(0.5;a,b)]^{1/c}\right)\right\}^{-1/\theta}.
\end{equation}

	\subsection{Moments and Moment Generating Function}

\begin{figure}[ht]
\centering
\subfloat[]
{\includegraphics[scale=0.4]
{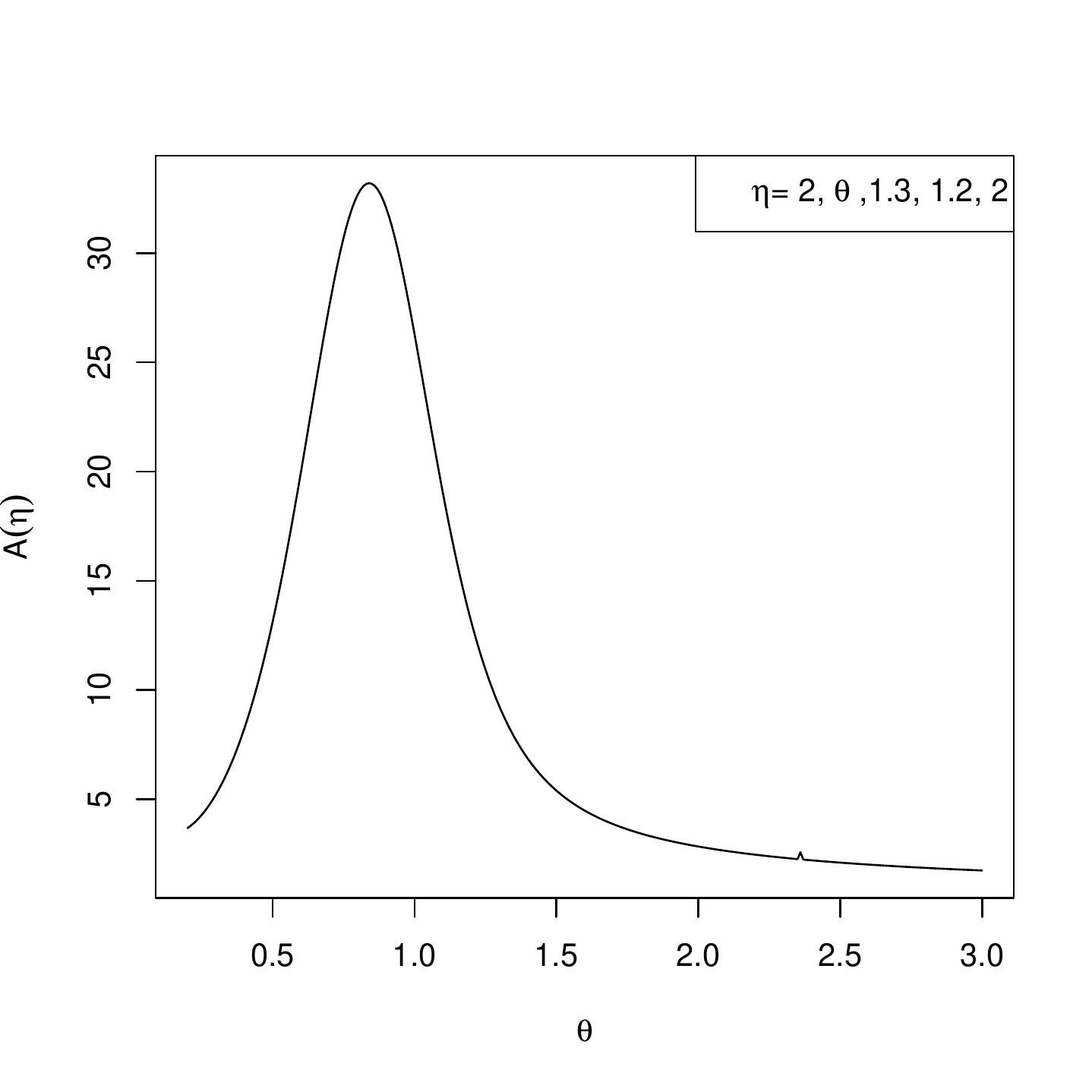} \label{asimm_theta}}
\subfloat[]
{\includegraphics[scale=0.4]
{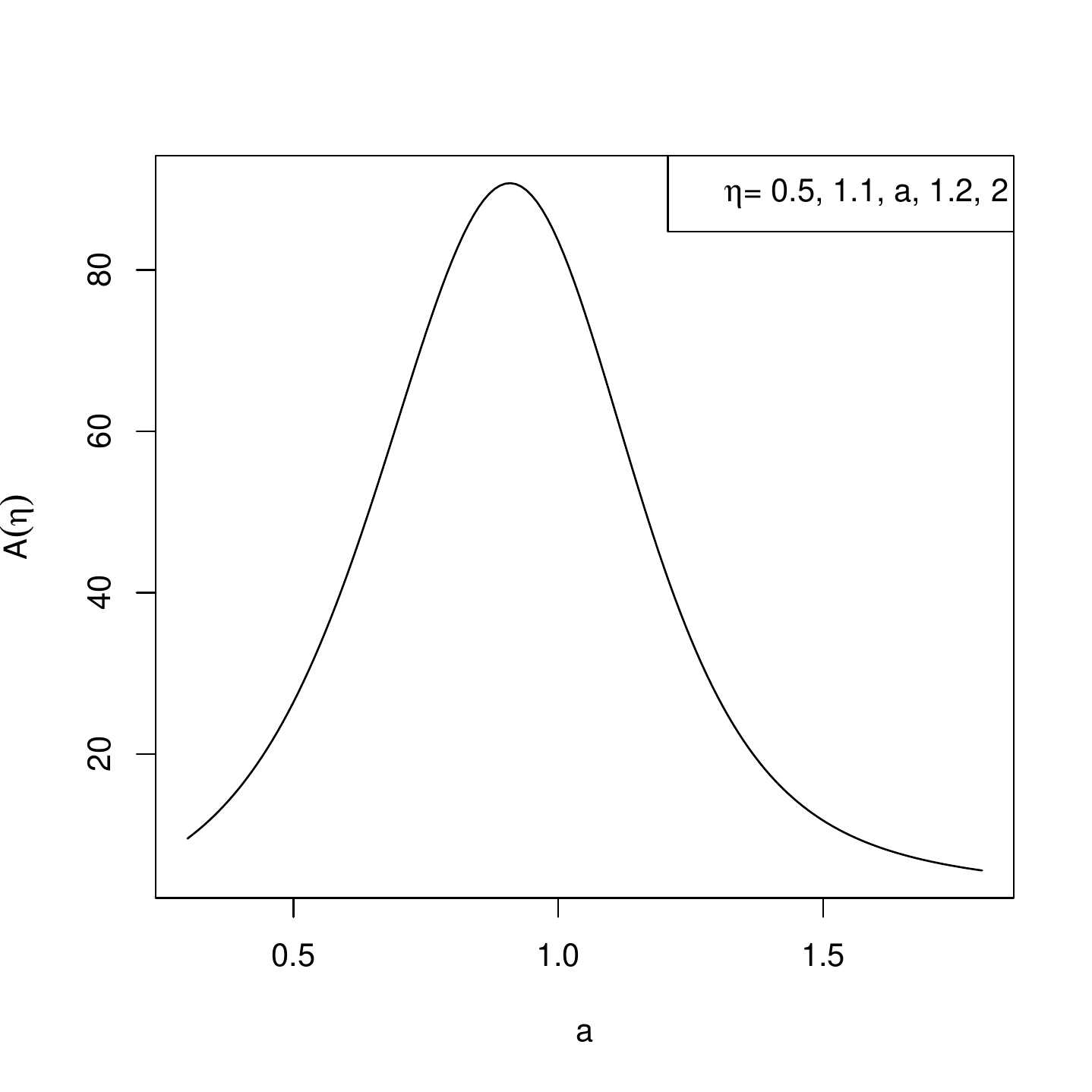}\label{asimm_a}} \\
\subfloat[]
{\includegraphics[scale=0.4]
{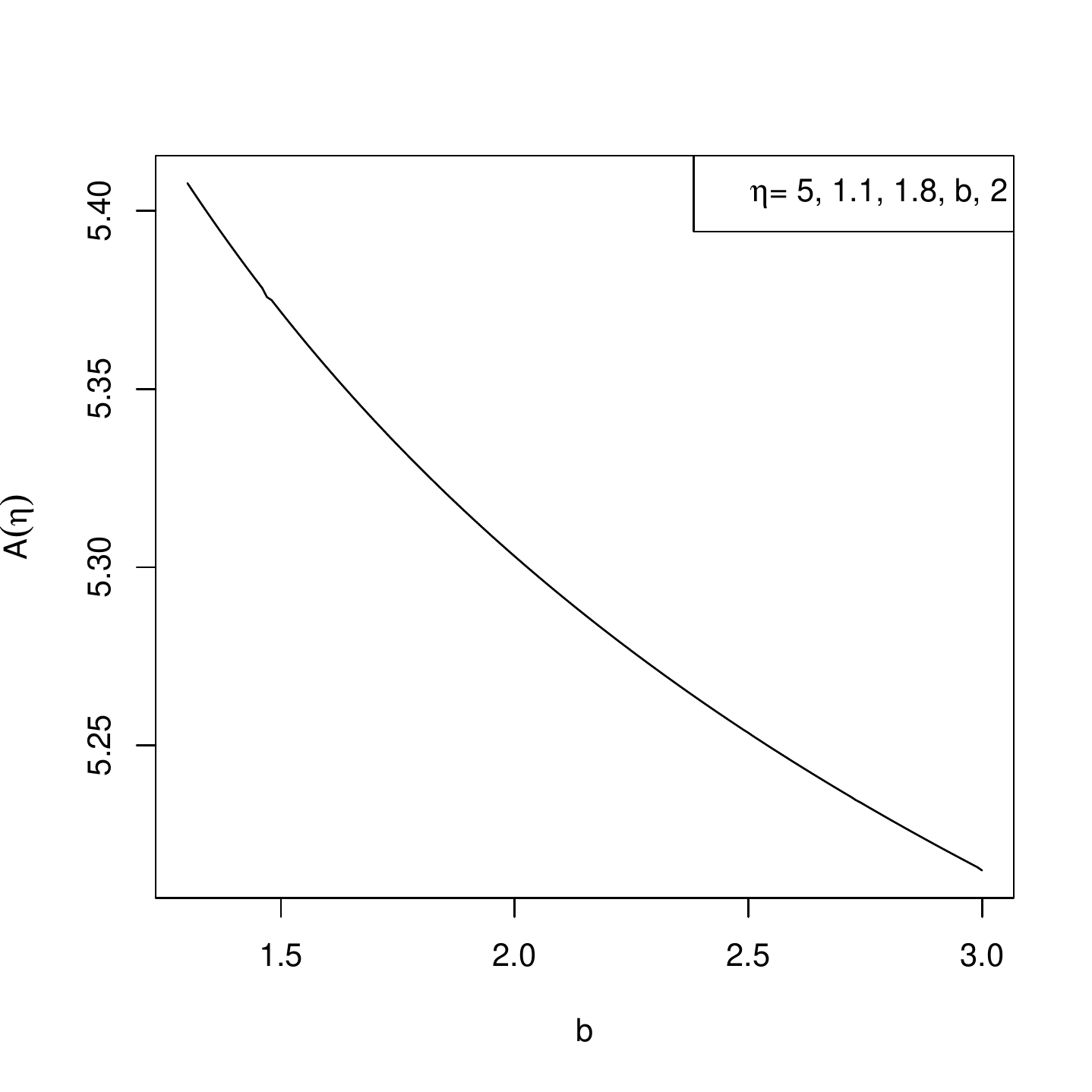}\label{asimm_b}}
\subfloat[]
{\includegraphics[scale=0.4]
{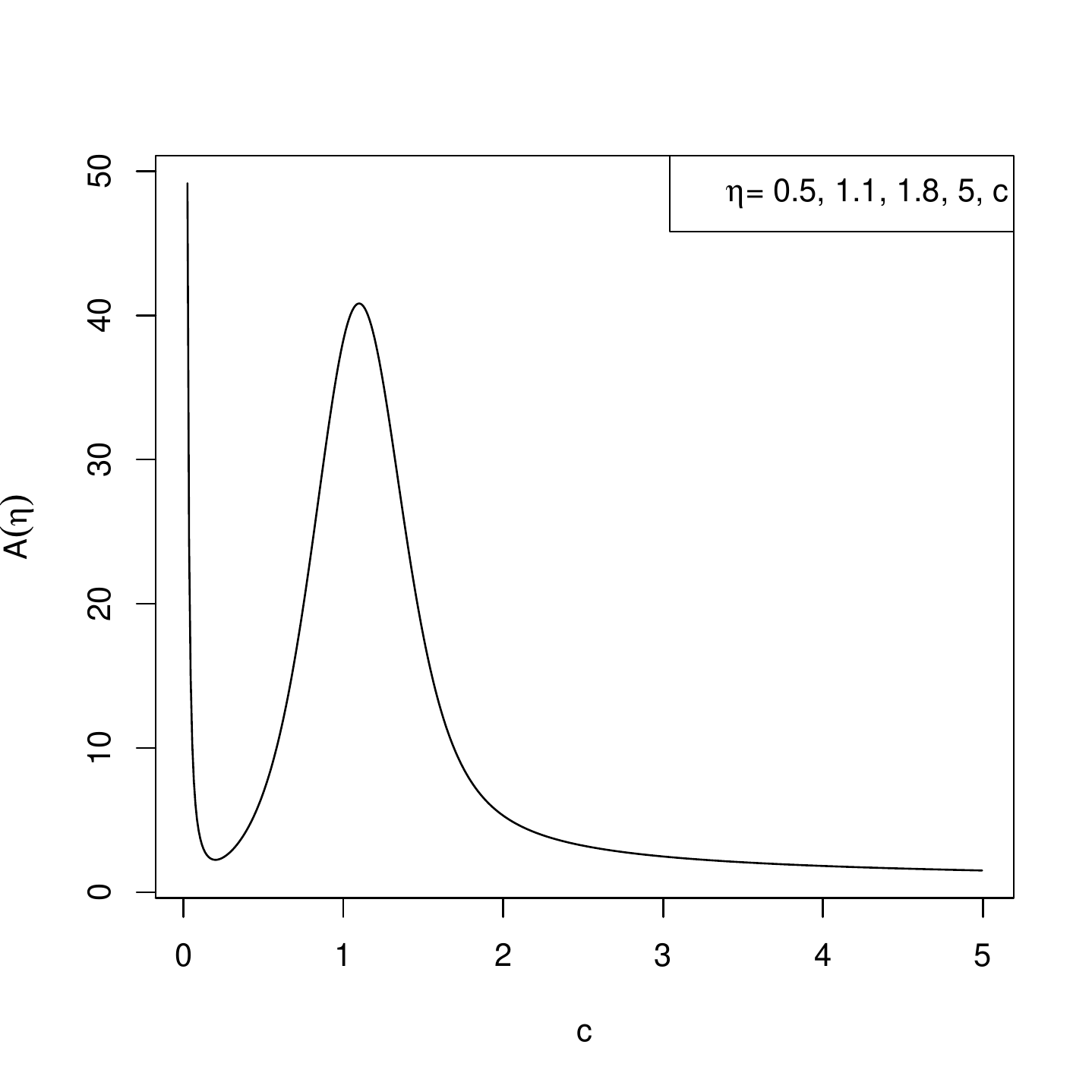}\label{asimm_c}}
\caption{Skewness measures of the \textit{rGB1-IWei} distribution for certain values of parameters.}
\label{fig:asimm}
\end{figure}

\begin{figure}[ht]
\centering
\subfloat[]
{\includegraphics[scale=0.4]
{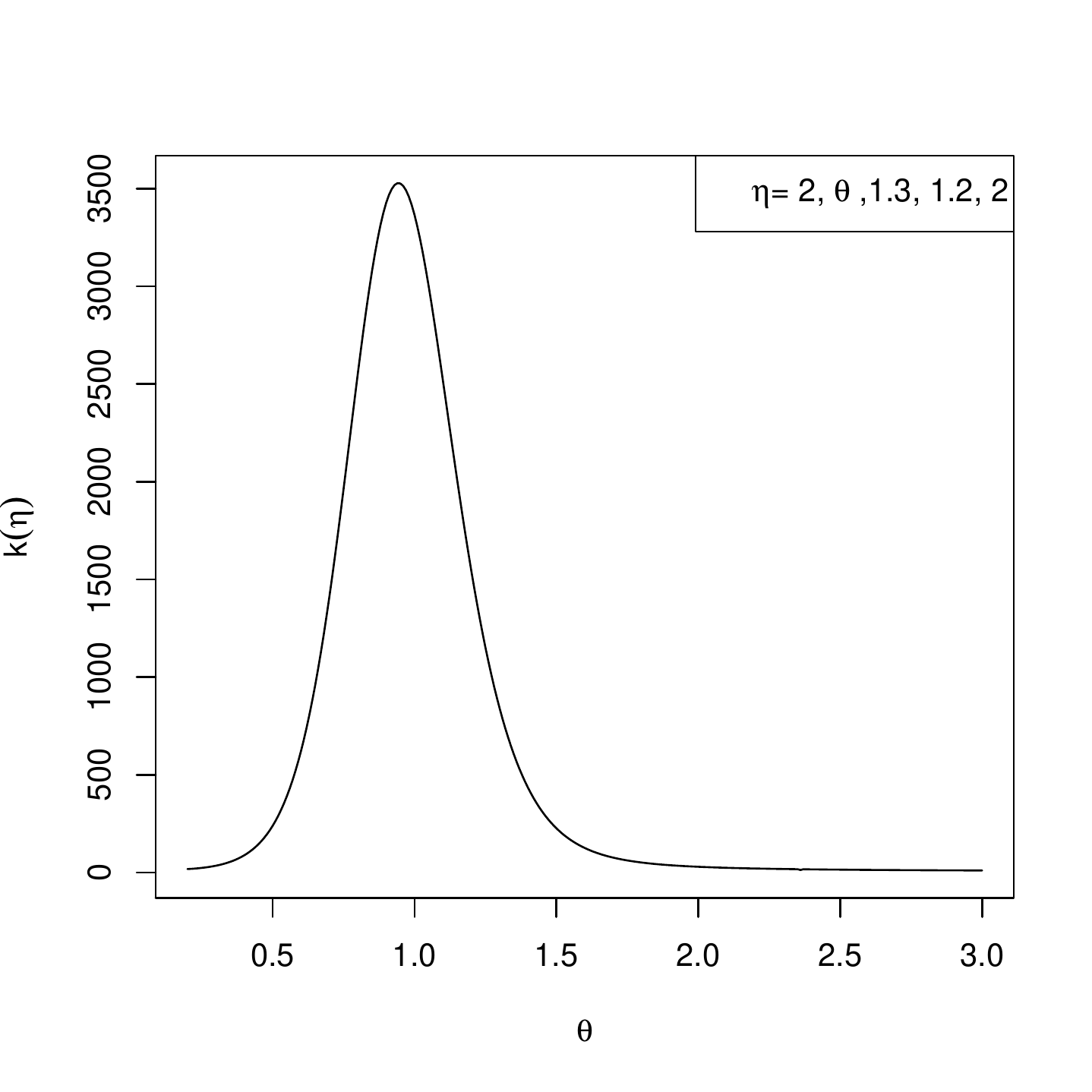} \label{kurt_theta}}
\subfloat[]
{\includegraphics[scale=0.4]
{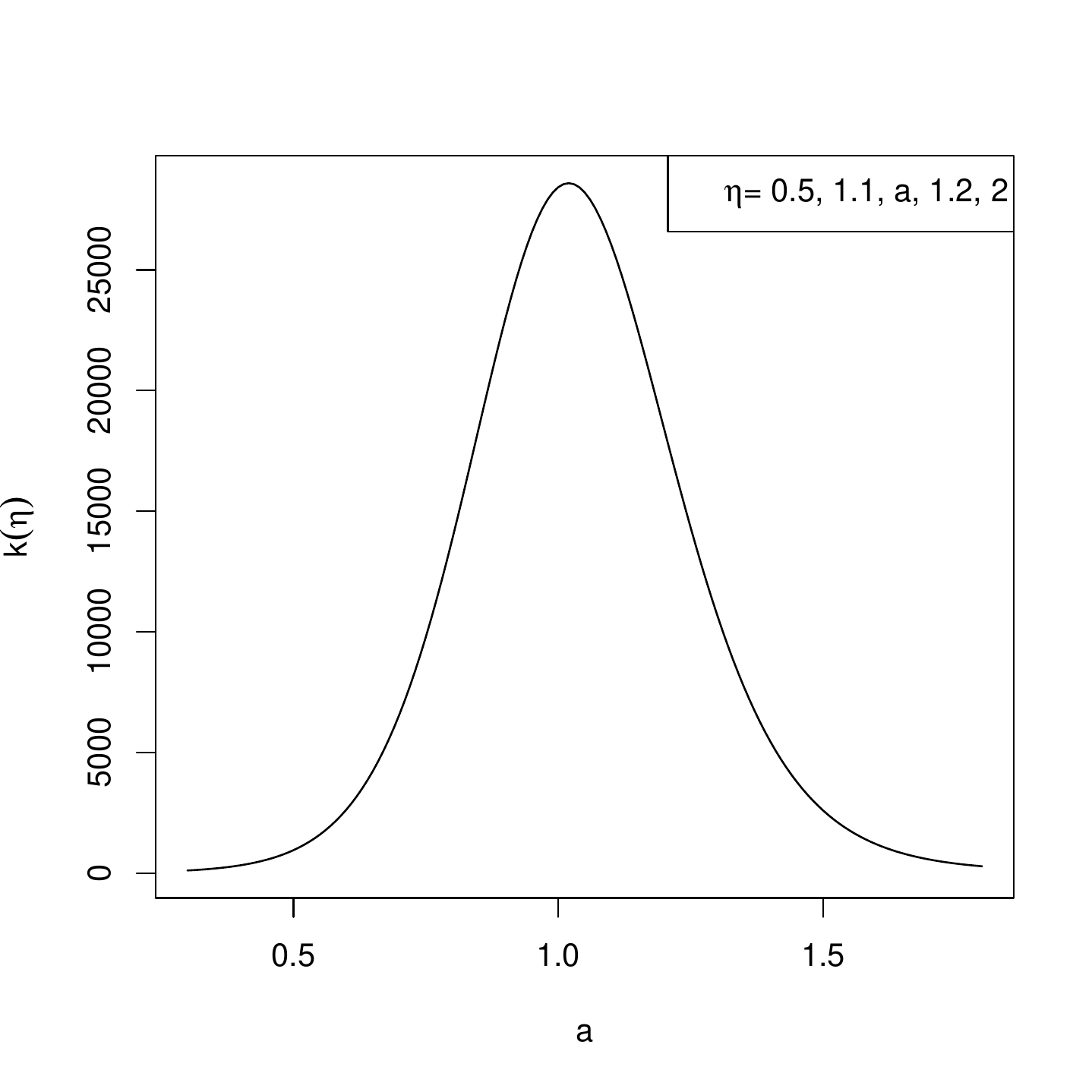}\label{kurt_a}} \\
\subfloat[]
{\includegraphics[scale=0.4]
{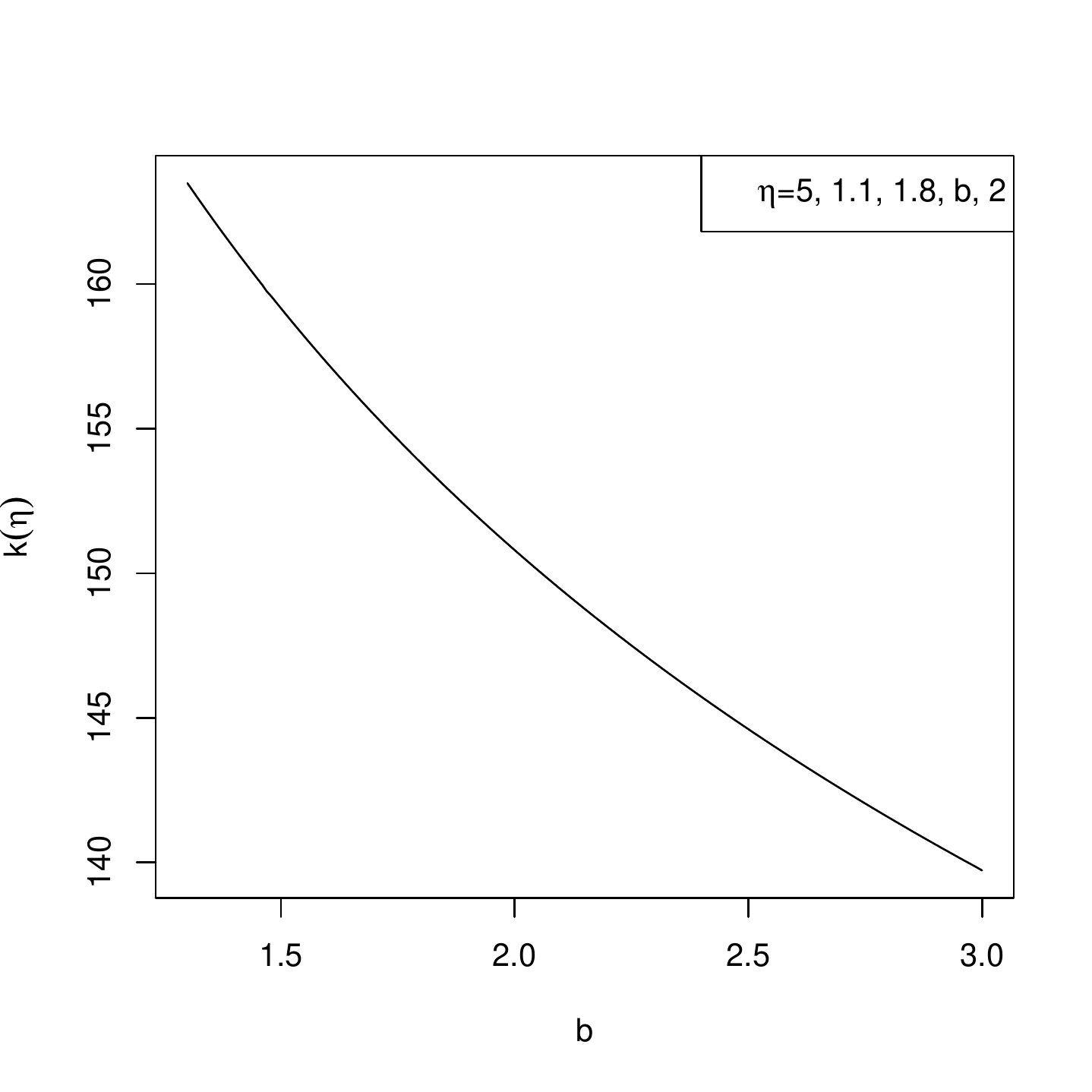}\label{kurt_b}}
\subfloat[]
{\includegraphics[scale=0.4]
{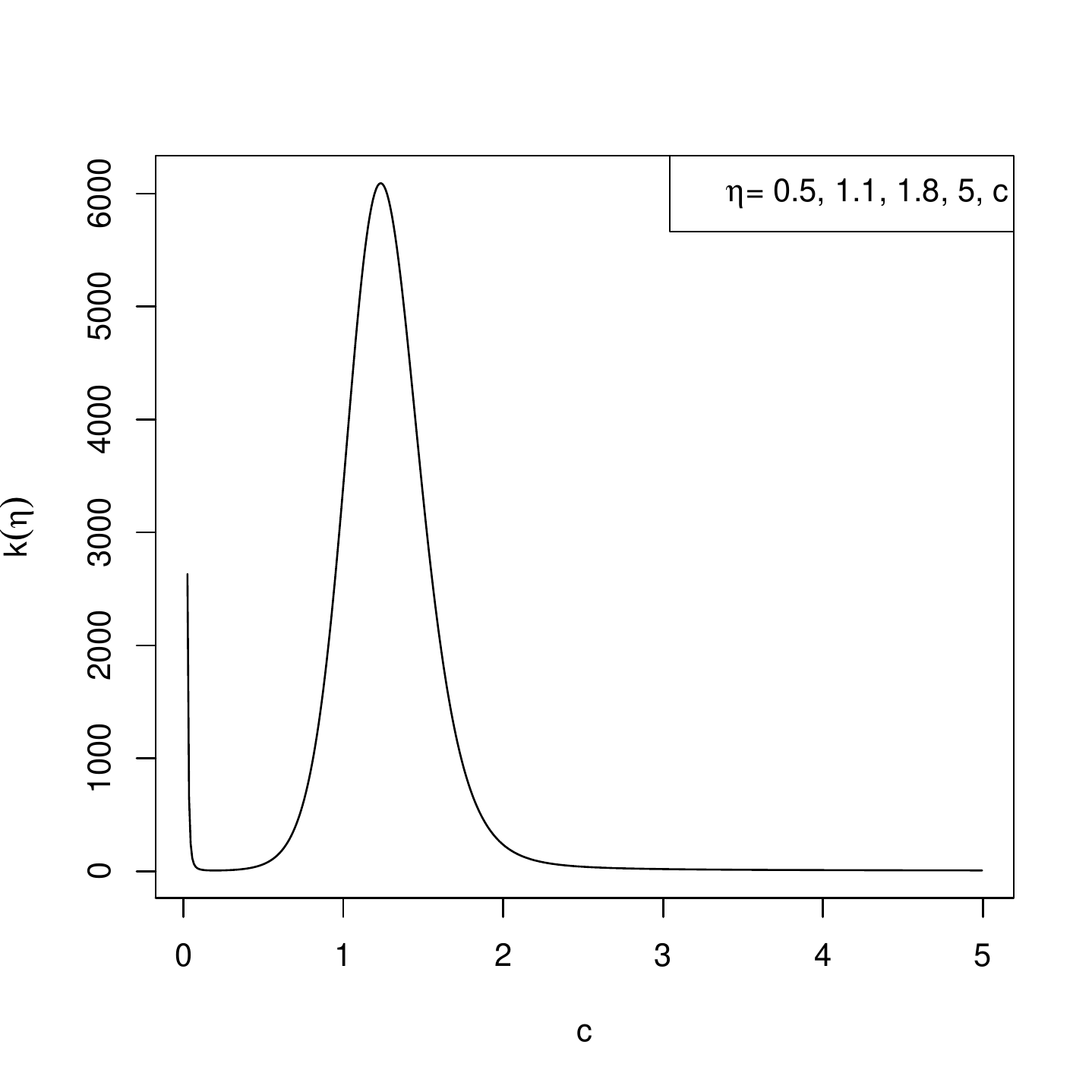}\label{kurt_c}}
\caption{Kurtosis measures of the \textit{rGB1-IWei} distribution for certain values of parameters.}
\label{fig:kur}
\end{figure}

In this subsection we discuss the $r_{th}$ moment for the \textit{rGB1-IWei} distribution. Moments are necessary in any statistical
analysis, especially in applications. They can be used to study the most important features and characteristics of a distribution (e.g., tendency,
dispersion, skewness and kurtosis).\medskip \newline
\begin{Prop}. If $X$ has \textit{rGB1-IWei} distribution, then the $r_{th}$ moment of $X$ is given by the following%
\begin{equation}
\mu _{r}(\boldsymbol{\eta})=\frac{c}{B(a,b)}\sum_{j_{1}=0}^{\infty
}p_{j_{1},b}\left( \sum_{j_{2}=0}^{\infty }\frac{p_{j_{2},c(a+j_{1})}}{%
(j_{2}+1)} \cdot \left[(j_{2}+1) \gamma \right] ^  {\frac{r}{\theta }} \cdot  \Gamma \left(1-\frac{r}{\theta }\right)\right)
\label{mo}
\end{equation}%
\end{Prop}
\textbf{\textit{Proof}}. It is an immediate consequence of the expansions of the \textit{pdf} defined in (\ref{pdfrGB1-GIWei_2}) and by equation (\ref{M_GIWei}). \\
Based on the first four moments the measures of skewness $A(\boldsymbol{\eta} )$ and kurtosis $k(\boldsymbol{\eta} ) $ of the \textit{rGB1-IWei} distribution can obtained
as%
\begin{equation*}
A(\boldsymbol{\eta} )=\frac{\mu _{3}(\boldsymbol{\eta} )-3\mu _{1}(\boldsymbol{\eta} )\mu _{2}(\boldsymbol{\eta} )+2\mu
_{1}^{3}(\boldsymbol{\eta} )}{\left[ \mu _{2}(\boldsymbol{\eta} )-\mu _{1}^{2}(\boldsymbol{\eta} )\right] ^{
\frac{3}{2}}},
\end{equation*}%
and
\begin{equation*}
k(\boldsymbol{\eta} )=\frac{\mu _{4}(\boldsymbol{\eta} )-4\mu _{1}(\boldsymbol{\eta} )\mu _{3}(\boldsymbol{\eta} )+6\mu
_{1}^{2}(\boldsymbol{\eta} )\mu _{2}(\boldsymbol{\eta} )-3\mu _{1}^{4}(\boldsymbol{\eta} )}{\left[ \mu
_{2}(\boldsymbol{\eta} )-\mu _{1}^{2}(\boldsymbol{\eta} )\right] ^{2}}.
\end{equation*}
Plots of the skewness and kurtosis measures of the \textit{rGB1-IWei} distribution for selected values of parameters are reported in Fig. \ref{fig:asimm} and Fig. \ref{fig:kur}.

Now, we shall derive the moment generating function of \textit{rGB1-IWei} distribution. \\
\begin{Prop}. If $X$ has \textit{rGB1-IWei} distribution, then the moment generating function $M_{X}(t)$ has the following form%
\begin{equation}
M_{X}(t)=\sum_{r=0}^{\infty }\left( \frac{t^{r}}{r!}\right) \frac{c}{%
B(a,b)}\sum_{j_{1}=0}^{\infty }p_{j_{1},b}\left( \sum_{j_{2}=0}^{\infty }%
\frac{p_{j_{2},c(a+j_{1})}}{(j_{2}+1)}\right) \left[\gamma (j_{2}+1) \right] ^  {\frac{r}{\theta }}
 \Gamma \left(1-\frac{r}{\theta }\right).
\label{fgm}
\end{equation}%
\end{Prop}
\textbf{\textit{Proof}}. Using the series expansion $e^{tx}=\sum_{r=0}^{\infty }\frac{(tx)^{r}}{r!}$, we can write
\begin{eqnarray}
M_{X}(t) &=&\sum_{r=0}^{\infty }\frac{(t)^{r}}{r!}
\int_0^{\infty }x^{r}f_{rGB1-IWei}(x;\boldsymbol{\eta })dx = \sum_{r=0}^{\infty }\frac{(t)^{r}}{r!}\mu _{r}(\boldsymbol{\eta }).
\label{fgm1}
\end{eqnarray}
 Substituting (\ref{mo}) into (\ref{fgm1}) we get the expression of $M_{X}(t)$. 


	\subsection{Entropy}
	
	The concept of entropy plays a vital role in information theory. The entropy of a random variable is defined in terms of its probability distribution and
can be shown to be a valid measure of randomness or uncertainty. In this section we present R\'{e}nyi entropy $I_{R}(\rho )$ \ and Shannon entropy $%
H(f_{rGB1-IWei}(x;\boldsymbol{\eta }))$.\medskip \newline

\begin{Prop}. If $X$ has \textit{rGB1-IWei} distribution, then the R\'{e}nyi entropy is given by
\begin{eqnarray}
I_{R}(\rho )&=&\frac{\rho}{1-\rho } \log \left(\frac{c}{B(a,b)}\right) + \frac{1}{1-\rho } \log \Gamma\left(\rho + \frac{(\rho-1)}{\theta } \right) -\log (\gamma \theta) \nonumber \\
&+& \frac{1}{1-\rho } \log \left( \sum^{\infty}_{j_{1}=0} \sum ^{\infty}_{j_{2}=0} \frac{p_{j_{1}, \rho(b-1)+1} p_{j_{2}, \rho(ac-1)+j_{1}c+1}    }
{ \left[ (j_{2}+\rho )\gamma \right] ^{ \rho+\frac{(\rho-1)}{\theta }  }  }
\right)
\label{IR}
\end{eqnarray}
\end{Prop}
\textbf{\textit{Proof}}. The R\'{e}nyi entropy is defined as%
\begin{equation*}
I_{R}(\rho )=\frac{1}{1-\rho }\log \left[ \int_{0}^{\infty }f_{rGB1-IWei}(x;\boldsymbol{\eta })^{\rho }dx\right] ,
\end{equation*}%
where $\rho >0,$and $\rho \neq 1$. For the \textit{rGB1-IWei} distribution, using the power series representation, the integral in $I_{R}(\rho )$ can be reduced to%
\begin{eqnarray*}
&&\int_{0}^{\infty }f_{rGB1-IWei}(x;\boldsymbol{\eta })^{\rho }dx = \\
&&=\left(\frac{c\gamma\theta}{B(a,b)}\right)^{\rho} \sum^{\infty}_{j_{1}=0} \sum ^{\infty}_{j_{2}=0} p_{j_{1}, \rho(b-1)+1} p_{j_{2}, \rho(ac-1)+j_{1}c+1}
\int^{\infty}_{0} x^{-\rho(\theta+1)} e^{-(j_{2}+\rho)\gamma x^{-\theta}}dx.
\end{eqnarray*}%
Setting $\gamma (j_{2}+\rho) x^{-\theta }=w$, after algebra, we obtain%
\begin{eqnarray*}
\int_{0}^{\infty }f_{rGB1-IWei}(x;\boldsymbol{\eta })^{\rho }dx
=\left(\frac{c\gamma\theta}{B(a,b)}\right)^{\rho} \frac{\Gamma\left(\rho+ \frac{(\rho-1)}{\theta}  \right)}{\gamma \theta}
\sum^{\infty}_{j_{1}=0} \sum ^{\infty}_{j_{2}=0} \frac{p_{j_{1}, \rho(b-1)+1} p_{j_{2}, \rho(ac-1)+j_{1}c+1}}{\left[(j_{2}+\rho )\gamma\right]^{\rho+\frac{(\rho-1)}{\theta}} }
\end{eqnarray*}%
through simple algebra we obtain (\ref{IR}).

\begin{Prop}. If $X$ has \textit{rGB1-IWei} distribution, then the Shannon entropy is given by
\begin{eqnarray}
H(f_{rGB1-IWei}(x;\boldsymbol{\eta }))&=& -\log\left(\frac{c \theta \gamma}{B(a,b)}\right)+\frac{c}{B(a,b)} \sum^{\infty}_{j_{1}=0} \sum ^{\infty}_{j_{2}=0} p_{j_{1}, b} p_{j_{2}, c(a+j_{1})} \nonumber \\
&& \left\{ \frac{(\theta+1)}{\theta} \left( \frac{\log\left[ (j_{2}+1)\gamma \right] - \Gamma{'}(1)}{(j_{2}+1)}\right)-
\frac{(ac-1)}{(j_{2}+1)} \left[\psi(1)-\psi(j_{2}+2)\right] + \right. \nonumber  \\
&& \left. +\frac{1}{(j_{2}+1)^2}-(b-1) \sum^{\infty}_{i=0} \frac{p_{i,j_{2}+1}}{(i+1)}
\left[ \psi(1)-\psi \left( \frac{i+1}{c}+1 \right)   \right]
\right\}
\end{eqnarray}
\end{Prop}

\textbf{\textit{Proof}}. The Shannon entropy is given by%
\begin{eqnarray*}
H(f_{rGB1-IWei}(x;\boldsymbol{\eta })) &=& - E \left[ \log f_{rGB1-IWei}(x;\boldsymbol{\eta }) \right]   \\
&=&-\log (c \gamma\theta  )+\log B(a,b)+(\theta +1)E \left[\log (X)\right] \\
&&+E\left[ \gamma X^{-\theta }\right] -(ac-1)E\left\{ \log \left[ 1-e^{-\gamma  X^{-\theta }}\right] \right\}  \\
&&-(b-1)E\left\{ \log \left( 1-\left[ 1-e^{-\gamma X^{-\theta }}\right]^{c} \right)  \right\}.
\end{eqnarray*}%
Using the power series representation, we have%
\begin{eqnarray*}
E\left[\log (X)\right] &=& \frac{c \theta \gamma}{B(a,b)} \sum^{\infty}_{j_{1}=0} \sum ^{\infty}_{j_{2}=0} p_{j_{1}, b} p_{j_{2}, c(a+j_{1})}
\int_{0}^{\infty }\log (x) x^{-(\theta+1)} e^{-(j_{2}+1)\gamma x^{-\theta} } dx \\
&=& \frac{c}{B(a,b)} \sum^{\infty}_{j_{1}=0} \sum ^{\infty}_{j_{2}=0} \frac{ p_{j_{1}, b} p_{j_{2}, c(a+j_{1})} }{\theta ( j_{2}+1)}
\int_{0}^{\infty } \left\{ \log\left[ ( j_{2}+1)\gamma \right] - \log(w)  \right\} e^{-w} dw =\\
&=& \frac{c}{B(a,b)} \sum^{\infty}_{j_{1}=0} \sum ^{\infty}_{j_{2}=0} \frac{ p_{j_{1}, b} p_{j_{2}, c(a+j_{1})} }{\theta( j_{2}+1)}
\left\{ \log\left[ ( j_{2}+1 )\gamma \right] - \Gamma^{'}(1)    \right\}
\end{eqnarray*}%
also, by \textit{Proposition 1}, we can write%
\begin{eqnarray*}
E\left[ \gamma X^{-\theta }\right] &=& \gamma E\left(X ^{-\theta} \right) =
\frac{c}{B(a,b)} \sum^{\infty}_{j_{1}=0} \sum ^{\infty}_{j_{2}=0} \frac{p_{j_{1}, b} p_{j_{2}, c(a+j_{1})}}{ ( j_{2}+1)^2  }.
\end{eqnarray*}%

Moreover, we have

\begin{eqnarray*}
E\left\{ \log \left[ 1-e^{-\gamma X^{-\theta }}\right] \right\} &=&
\frac{c \theta \gamma}{B(a,b)} \sum^{\infty}_{j_{1}=0} \sum ^{\infty}_{j_{2}=0} p_{j_{1}, b} p_{j_{2}, c(a+j_{1})} \\
&& \int_{0}^{\infty } \log \left[ 1-e^{-\gamma x^{-\theta }}\right]  x^{-(\theta+1)} \left[ e^{-\gamma x^{-\theta} }  \right] ^{j_{2}+1} dx,
\end{eqnarray*}%
setting $w=1-e^{-\gamma x^{-\theta }}$, after algebra, we obtain
\begin{eqnarray*}
E\left\{ \log \left[ 1-e^{-\gamma x^{-\theta }}\right] \right\} &=&
\frac{c}{B(a,b)} \sum^{\infty}_{j_{1}=0} \sum ^{\infty}_{j_{2}=0} p_{j_{1}, b} p_{j_{2}, c(a+j_{1})} \int^{1}_{0} \log(w) (1-w)^{j_{2}} dw =    \\
&=& \frac{c}{B(a,b)} \sum^{\infty}_{j_{1}=0} \sum ^{\infty}_{j_{2}=0} \frac{p_{j_{1}, b} p_{j_{2}, c(a+j_{1})}}{ (j_{2}+1) }
\left\{ \psi(1)-\psi(2+j_{2})  \right\},
\end{eqnarray*}
where $\psi(.)$ denote the digamma function. Finally, the last integral is given by
\begin{eqnarray*}
&& E\left\{ \log  \left( 1-  \left[ 1-e^{-\gamma X^{-\theta }}\right]^c \right)       \right\} = \\
&=& \frac{c \theta \gamma}{B(a,b)} \sum^{\infty}_{j_{1}=0} \sum ^{\infty}_{j_{2}=0} p_{j_{1}, b} p_{j_{2}, c(a+j_{1})}
\int^{\infty}_{0} \log  \left( 1-  \left[ 1-e^{-\gamma x^{-\theta }}\right]^c \right) x^{-(\theta+1)} \left[ e^{-\gamma x^{-\theta} }  \right] ^{j_{2}+1} dx= \\
&=& \frac{1}{B(a,b)} \sum^{\infty}_{j_{1}=0} \sum ^{\infty}_{j_{2}=0} p_{j_{1}, b} p_{j_{2}, c(a+j_{1})}
\int^{1}_{0} \log(w) (1-w)^{\frac{1}{c}-1 } \left[ 1- (1-w)^{\frac{1}{c} } \right]^{j_{2}} dw = \\
&=& \frac{1}{B(a,b)} \sum^{\infty}_{j_{1}=0} \sum ^{\infty}_{j_{2}=0} p_{j_{1}, b} p_{j_{2}, c(a+j_{1})}
\sum^{\infty}_{i=0} p_{i, j_{2}+1} \int^{1}_{0} \log(w) (1-w)^{ \frac{i+1}{c} -1  } dw=  \\
&=&  \frac{c}{B(a,b)} \sum^{\infty}_{j_{1}=0} \sum ^{\infty}_{j_{2}=0} p_{j_{1}, b} p_{j_{2}, c(a+j_{1})}
\sum^{\infty}_{i=0} \frac{p_{i, j_{2}+1}}{(i+1)} \left\{ \psi(1)-\psi \left( \frac{i+1}{c} +1 \right)  \right\}.
\end{eqnarray*}

	\subsection{Order Statistics}
	
In this subsection, we derive closed form expressions for the \textit{pdf} and moments of the $k_{th}$ order statistic of the \textit{rGB1-IWei} distribution. Let $X_{1},X_{2},...,X_{n}$ be a simple random sample from \textit{rGB1-IWei} distribution with \textit{cdf} and \textit{pdf}  given by (\ref{dfrGB1-GIWei}) and (\ref{pdfrGB1-GIWei}), respectively.
Let $X_{1:n},X_{2:n},...,X_{n:n}$ denote the order statistics obtained from this sample. We now give the probability density function of $X_{k:n}$, say $%
f_{_{k:n}rGB1-IWei}(x;\boldsymbol{\eta })$ and the moments of $X_{k:n}$, $k=1,2,...,n$. The probability density function of $X_{k:n}$ is given by%
\begin{eqnarray}
f_{_{k:n}rGB1-IWei}(x;\boldsymbol{\eta })&=&\frac{1}{B(k,n-k+1)}\left[
F_{rGB1-IWei}(x;\boldsymbol{\eta })\right] ^{k-1} \nonumber \\
&& \left[ 1-F_{rGB1-IWei}(x;\boldsymbol{\eta }))\right] ^{n-k} f_{rGB1-IWei}(x;\boldsymbol{\eta })
\end{eqnarray}%
where $F_{rGB1-IWei}(x;\boldsymbol{\eta })$ and $f_{rGB1-IWei}(x;\boldsymbol{\eta })$ are the \textit{cdf} and \textit{pdf} of the \textit{rGB1-IWei} distribution given by (\ref{dfrGB1-GIWei}) and (\ref{pdfrGB1-GIWei}), respectively. Since $0<F_{rGB1-IWei}(x;\boldsymbol{\eta })<1 $, for $x>0$, by using the binomial series expansion of $\left[1-F_{rGB1-IWei}(x;\boldsymbol{\eta })\right] ^{n-k}$, given by%
\begin{equation*}
\left[ 1-F_{rGB1-IWei}(x;\boldsymbol{\eta })\right] ^{n-k}=\sum_{j=0}^{n-k}(-1)^{j}\binom{n-k}{j}\left[ F_{rGB1-IWei}(x;\boldsymbol{%
\eta })\right] ^{^{j}},
\end{equation*}%
we have%
\begin{equation}
f_{_{k:n}rGB1-IWei}(x;\boldsymbol{\eta })=\sum_{j=0}^{n-k}(-1)^{j}%
\binom{n-k}{j}\left[ F_{rGB1-IWei}(x;\boldsymbol{\eta })\right]
^{k+j-1}f_{rGB1-IWei}(x;\boldsymbol{\eta }),
\label{pdfOS}
\end{equation}%
substituting (\ref{dfrGB1-GIWei}) and (\ref{pdfrGB1-GIWei}) into (\ref{pdfOS}), we can express the $s_{th}$ ordinary moment of the $k_{th}$ order statistics $X_{k:n}$, say  $E(X_{k:n}^{s})$, as a liner combination of the probability weighted moments of the \textit{rGB1-IWei} distribution, \textit{i.e.}

\begin{equation}
E \left[X_{k:n}^{s} \right] = \sum_{j=0}^{n-k}(-1)^{j} \binom{n-k}{j} E_{rGB1-IWei}  \left[ X^{s} \left( F_{rGB1-IWei}(X;\boldsymbol{\eta })\right)^{k+j-1} \right].
\label{MomOS}
\end{equation}%

	\subsection{Reliability}
	
In reliability studies, the stress-strength term is often used to describe the life of a component which has a random strength $X$ and is subject to a random stress $Y$. The component fails if the stress applied to it exceeds the strength, and the component will function satisfactorily whenever $Y<X$. Thus $R=P(Y<X)$ is a measure of a component reliability which has many applications in physics, engineering, genetics, psychology and economics (Kotz et al. 2003). In this subsection, we derive the reliability $R$ when $X\sim rGB1-IWei(\boldsymbol{\eta }_x)$ and $Y\sim rGB1-IWei(\boldsymbol{\eta }_y)$ are independent random variables. \\
Using expansion (\ref{pdfrGB1-GIWei_2}) of the \textit{pdf} of the \textit{rGB1-IWei} random variable, after algebra, we can write
\begin{eqnarray*}
R=\frac{c_{x}c_{y} }{B(a_{x},b_{x})B(a_{y},b_{y}) } \sum^{\infty }_{j_{1}=0} \sum^{\infty }_{j_{2}=0} \sum^{\infty }_{j_{3}=0} \sum^{\infty }_{j_{4}=0}
\frac{ p_{j_{1},b_y} p_{j_{2},c_{y} ( a_{y} +j_{1})} p_{j_{3},b_x} p_{j_{4},c_{x} ( a_{x} +j_{3})} }
{ ( j_{2}+1) ( j_{4}+1) }  R_{j_{2},j_{4} }
\end{eqnarray*}
where
\begin{eqnarray}
R_{j_{2},j_{4} }=\int^{\infty}_{0} g_{IWei}(x; \gamma _{x}(j_{4}+1), \theta_{x} ) G_{IWei}(x; \gamma _{y}(j_{2}+1),  \theta_{y} ) dx
\end{eqnarray}
is the reliability between the independent Inverse Weibull random variables. Hence, the reliability between \textit{rGB1-IWei} random variables is a linear combination of the reliabilities between \textit{IWei} random variables. \\
Putting $\theta_{x}=\theta_{y}=\theta $, after algebra, the $ R_{j_{2},j_{4} } $ is given by
\begin{eqnarray*}
R_{j_{2},j_{4} }=\gamma_{x}(j_{4}+1) \theta  \int^{\infty}_{0} x^{-\theta-1} e^{-x^{-\theta}
\left[ \gamma_{x}(j_{4}+1)  + \gamma_{y}(j_{2}+1)   \right]   } dx 
= \frac{\gamma_{x}(j_{4}+1)  }
{\left[ \gamma_{x}(j_{4}+1)  + \gamma_{y}(j_{2}+1)   \right]}.
\end{eqnarray*}
Finally, the reliability between independent \textit{rGB1-IWei} random variables, in the case $\theta_{x}=\theta_{y}=\theta $, is given by
\begin{eqnarray*}
R=\frac{c_{x}c_{y} }{B(a_{x},b_{x})B(a_{y},b_{y}) } \sum^{\infty }_{j_{1}=0} \sum^{\infty }_{j_{2}=0} \sum^{\infty }_{j_{3}=0} \sum^{\infty }_{j_{4}=0}
\frac{    p_{j_{1},b_y} p_{j_{2},c_{y} ( a_{y} +j_{1})} p_{j_{3},b_x} p_{j_{4},c_{x} ( a_{x} +j_{3})} \gamma_{x}(j_{4}+1)  }
{ ( j_{2}+1) ( j_{4}+1) \left[ \gamma_{x}(j_{4}+1)  + \gamma_{y}(j_{2}+1)  \right]   } .
\end{eqnarray*}

\section{Link between Generalizations of Inverse Weibull and  Generalizations of Dagum distribution }

In this section, we generalize the link between the Inverse Weibull distribution and Dagum distribution, via mixing Gamma density (see, for example,\citet[pag.192]{Kleiber2003}). \\
In the literature, it is recognized that a \textit{pdf}, say $p_1(x)$, has a mixture representation if it can be written as
$p_1(x)=\int_{\Theta} p_2(x|\theta)p_3(\theta) d\theta$, where $\theta$ is regarded as a random variable with \textit{pdf} $p_3(.)$, called mixing density, and $p_2(x|\theta)$ is the conditional \textit{pdf} of \textit{rv} $X$ given $\theta$. Mixtures are also called \textit{compound distributions}. The following \textit{Proposition} shows that the \textit{Burr III} distribution (or special case of the Dagum distribution) has a mixture representation with mixing Gamma density.

\begin{Prop}. If the conditional \textit{rv} $X$ given $\gamma$ has Inverse Weibull distribution, \textit{i.e.} $X|\gamma \sim IWei(\gamma, \theta) $, and
the \textit{rv} $\gamma$ is Gamma distributed, \textit{i.e.} $\gamma \sim Ga(\beta)$, then the \textit{rv} $X$ is distributed as a \textit{Burr III} distribution, \textit{i.e.} $X\sim BurrIII(\beta,\theta)$.
\end{Prop}
\textbf{\textit{Proof}}. \textit{See, for example, Kotz and Kleiber (2003), pag. 192.} \\


Next, by means of two propositions, we shall discuss the connection between certain generalizations of the \textit{IWei} distribution and those of the Dagum distribution via mixing Gamma density.

\begin{Prop}. If $X|\gamma \sim BeIWei(a,b,\gamma, \theta)$ and $\gamma \sim Ga(\beta)$ then $X\sim \mbox{Mixture of Dag}(\beta,\lambda_j,\theta)$, with $\lambda_j=(a+j)$.
\end{Prop}
\textbf{\textit{Proof}}. The conditional \textit{pdf} of the \textit{rv} $X|\gamma$ and the marginal \textit{pdf} of the \textit{rv} $\gamma$, respectively, are
\begin{equation*}
f(x|\gamma; a,b, \theta)=\frac{\gamma \theta }{B(a,b)} x^{-\theta-1} e^{-\gamma x^{-\theta} }
\left[e^{-\gamma x^{-\theta} }\right]^{a-1}
\left[  1- e^{-\gamma x^{-\theta} } \right]^{b-1}
\end{equation*}
and
\begin{equation*}
f(\gamma; \beta)=\frac{\gamma^{\beta -1} e^{-\gamma}}{\Gamma (\beta)}.
\end{equation*}
 Observed that
 $\left[1-e^{-\gamma x^{-\theta} } \right]^{b-1}= \sum^{\infty}_{j=0} \frac{(-1)^{j} \Gamma(b) }{\Gamma(b-j) j!} e^{-\gamma j x^{-\theta} } $, after algebra, the joint  \textit{pdf} between $X$ and $\gamma$ is
\begin{equation*}
f(x,\gamma; a,b,\beta,\theta)= \frac{ \theta x^{-\theta-1}} {B(a,b) \Gamma(\beta)}
\sum^{\infty}_{j=0} \frac{(-1)^{j} \Gamma(b) }{\Gamma(b-j) j!} \gamma^{\beta} e^{-\gamma \left[1+(a+j) x^{-\theta}  \right] }
\end{equation*}
marginalizing with respect to $\gamma$, we obtain the marginal \textit{pdf} of \textit{rv} $X$, \textit{i.e.}
\begin{eqnarray*}
f(x, a,b,\beta,\theta) &=& \frac{ \theta x^{-\theta-1}} {B(a,b) \Gamma(\beta)}
\sum^{\infty}_{j=0} \frac{(-1)^{j} \Gamma(b) }{\Gamma(b-j) j!} \int ^{\infty}_{0} \gamma^{\beta} e^{-\gamma \left[1+(a+j) (x^{-\theta}  \right] } d\gamma \\
&=& \sum^{\infty}_{j=0} \frac{(-1)^{j} \Gamma(a+b) }{\Gamma(a)\Gamma(b-j) j! (a+j)} f_{Dag} (x;\beta, (a+j),\theta)
\end{eqnarray*}

 where $f_{Dag} (x;\beta, (a+j),\theta)= \beta(a+j)\theta x^{-\theta-1} \left[1+(a+j) x^{-\theta} \right]^{-\beta-1}$ is the \textit{pdf} of a Dagum \textit{rv}. Moreover, given that $\sum^{\infty}_{j=0} \frac{(-1)^{j} \Gamma(a+b) }{\Gamma(a)\Gamma(b-j) j! (a+j)}=1$, see Domma and Condino (2013), we can say that $f(x, a,b,\beta,\theta)$ is a mixture of Dagum distribution.
 \hspace{12 cm} $\triangleleft$

\begin{Prop}. If $X|\gamma \sim rGB1-IWei(a,b,c,\gamma, \theta)$ and $\gamma \sim Ga(\beta)$ then $X\sim rGB1-Da(a,b,c,\beta, \lambda_j, \theta)$, with $\lambda_j=(a+j)$.
\end{Prop}
\textbf{\textit{Proof}}. The conditional \textit{pdf} of $X|\gamma$ is
\begin{eqnarray*}
f_{rGB1-IWei}(x|\gamma; a,b,c, \theta)=\frac{c}{B(a,b)}\sum ^{\infty}_{j_{1}=0} p_{j_{1},b}
\left(\sum ^{\infty}_{j_{2}=0}  \frac{p_{j_{2},c(a+j_{1})}}{(j_{2}+1)}   \gamma(j_{2}+1)\theta x^{-\theta-1}
e^ {-\gamma (j_{2}+1) x^{-\theta}}
\right).
\end{eqnarray*}
Assuming that the \textit{rv} $\gamma $ has \textit{pdf}  $f(\gamma; \beta)=\frac{\gamma^{\beta-1} e^{-\gamma}}{\Gamma(\beta)}$, the marginal \textit{pdf} of $X$ is given by
\begin{eqnarray*}
f(x;\boldsymbol{\eta}_1)&=&\int^{\infty}_{0}f_{rGB1-IWei}(x|\gamma; a,b,c,\alpha,\theta)f(\gamma; \beta) d\gamma \nonumber \\
&=& \frac{c}{B(a,b)}\sum ^{\infty}_{j_{1}=0} p_{j_{1},b}
\sum ^{\infty}_{j_{2}=0}  \frac{p_{j_{2},c(a+j_{1})}}{(j_{2}+1)}  \frac{\theta x^{-\theta-1}}{\Gamma(\beta)}
\int^{\infty}_{0} \gamma^{\beta}(j_{2}+1)
e^ {-\gamma \left[(j_{2}+1) x^{-\theta} +1\right]    } d\gamma \nonumber \\
&=& \frac{c}{B(a,b)}\sum ^{\infty}_{j_{1}=0} p_{j_{1},b}
\sum ^{\infty}_{j_{2}=0}  \frac{p_{j_{2},c(a+j_{1})}}{(j_{2}+1)}  \frac{\theta x^{-\theta-1}}{\Gamma(\beta)}
\frac{(j_{2}+1)}{\left[(j_{2}+1) x^{-\theta} +1\right] ^{\beta+1}} \int^{\infty}_{0}w^{\beta} e^{-w}dw \nonumber \\
&=& \frac{c}{B(a,b)}\sum ^{\infty}_{j_{1}=0} p_{j_{1},b}
\sum ^{\infty}_{j_{2}=0}  \frac{p_{j_{2},c(a+j_{1})}}{(j_{2}+1)} \beta  (j_{2}+1)  \theta x^{-\theta-1}
\left[(j_{2}+1)  x^{-\theta} +1\right] ^{-\beta-1} \nonumber \\
&=& \frac{c}{B(a,b)}\sum ^{\infty}_{j_{1}=0} p_{j_{1},b}
\sum ^{\infty}_{j_{2}=0}  \frac{p_{j_{2},c(a+j_{1})}}{(j_{2}+1)} f_{Da}(x;\beta, (j_{2}+1), \theta ).
\end{eqnarray*}
 \hspace{17 cm} $\triangleleft$

\section{Estimation and Inference}
In order to estimate the parameters $\boldsymbol{\eta}=(a,b,c,\gamma,\theta)$ of the \textit{rGB1-IWei} distribution, we use the maximum likelihood (\textit{ML}) method. Let $\textbf{x}=(x_1, x_2, ..., x_n)$ be a random sample of size $n$ from the \textit{rGB1-IWei} given by (\ref{pdfrGB1-GIWei}). The log-likelihood function for the vector of parameters $\boldsymbol{\eta}=(a,b,c,\gamma,\alpha,\theta)$ can be expressed as
\begin{eqnarray}
\ell(\boldsymbol{\eta})&=& n\log(c)-nlog(B(a,b))+(ac-1)\sum^{n}_{i=1} \log\left[1-G_{IWei}(x_i;\boldsymbol{\tau})\right] \nonumber \\
&+& (b-1) \sum^{n}_{i=1} \log\left\{ 1-\left[1-G_{IWei}(x_i;\boldsymbol{\tau})\right]^{c} \right\} +
		\sum^{n}_{i=1} \log\left\{ g_{IWei}(x_i;\boldsymbol{\tau})\right\}
\end{eqnarray}
where $\boldsymbol{\tau}=(\gamma, \theta)$. In what follows, we denote with $\dot{h}_{y} (y,z)=\frac{\partial h(y,z)}{\partial y}$, $\ddot{h}_{yy} (y,z)=\frac{\partial^2 h(y,z)}{\partial y^2}$ and $\ddot{h}_{yz} (y,z)=\frac{\partial^2 h(y,z)}{\partial y \partial z}$ the partial derivatives of first order, of second order and mixed of a function, say $h(y,z)$, respectively. Moreover, to simplify the notation, we use $G(x_i;\boldsymbol{\tau})=G_{IWei}(x_i;\boldsymbol{\tau})$ and $g(x_i;\boldsymbol{\tau})=g_{IWei}(x_i;\boldsymbol{\tau})$. \\
Differentiating $\ell(\boldsymbol{\eta})$ with respect to $a$, $b$, $c$, $\gamma$ and $\theta$, respectively, and setting the results equal to zero, we have
\begin{eqnarray*}
\left\{
\begin{array}{l}
\frac{\partial \ell (\boldsymbol{\eta} ) }{\partial a }=%
 -n \frac{ \dot{B}_{a} (a,b) }{B(a,b)} + c \sum ^{n}_{i=1} \lg \left[1-G(x_i;\boldsymbol{\tau})\right] =0
\\ \\

\frac{\partial \ell (\boldsymbol{\eta} )}{\partial b }=  %
-n \frac{ \dot{B}_{b} (a,b) }{B(a,b)} + \sum ^{n}_{i=1} \lg \left(1- \left[1-G(x_i;\boldsymbol{\tau})\right]^c \right) =0
 \\ \\

\frac{\partial \ell (\boldsymbol{\eta} )}{\partial c }=  %
\frac{n}{c}+ a \sum ^{n}_{i=1} \lg \left[1-G(x_i;\boldsymbol{\tau})\right] - (b-1) \sum ^{n}_{i=1}
\frac{\left[1-G(x_i;\boldsymbol{\tau})\right]^c \lg \left[1-G(x_i;\boldsymbol{\tau})\right]}
{ \left(1- \left[1-G(x_i;\boldsymbol{\tau})\right]^c \right) } =0
\\ \\

\frac{\partial \ell (\boldsymbol{\eta} )}{\partial \tau_{j} }=  %
- (ac-1) \sum ^{n}_{i=1} \frac{\dot{G}_{\tau_{j}} (x_i;\boldsymbol{\tau})}{\left[1-G(x_i;\boldsymbol{\tau})\right]}
+c(b-1) \sum ^{n}_{i=1}  \frac{\left[1-G(x_i;\boldsymbol{\tau})\right]^{c-1} \dot{G}_{\tau_{j}} (x_i;\boldsymbol{\tau}) }
{\left(1- \left[1-G(x_i;\boldsymbol{\tau})\right]^c \right)}+\sum ^{n}_{i=1} \frac{\dot{g}_{\tau_{j}}(x_i;\boldsymbol{\tau})}{g(x_i;\boldsymbol{\tau})}=0
 \end{array}
\right. \\
\end{eqnarray*}
for $j=1,2,3$ and with $\tau_1=\gamma$ and $\tau_2=\theta$, where $\dot{G}_{...} (x_i;\boldsymbol{\tau})$ and $\dot{g}_{...} (x_i;\boldsymbol{\tau})$ are reported in the \textit{Appendix}.
The system does not admit any explicit solutions; therefore, the \textit{ML} estimates $\boldsymbol{\hat{\eta}}=(\hat{a}, \hat{b}, \hat{c}, \hat{\gamma}, \hat{\theta})$ can only be obtained by means of numerical procedures. Under regularity conditions, the \textit{ML} estimator $\boldsymbol{\hat{\eta}}$ is consistent and asymptotically normally distributed. Moreover, the asymptotic variance-covariance matrix of $\boldsymbol{\hat{\eta}}$ can be approximated by the inverse of observed information matrix given by
\[
\boldsymbol {J}(\boldsymbol{\eta})=
  \begin{bmatrix}
    J_{aa} & J_{ab} & J_{ac} & J_{a\gamma}       & J_{a\theta} \\
    ..     & J_{bb} & J_{bc} & J_{b\gamma}       & J_{b\theta} \\
		..     & ..     & J_{cc} & J_{c\gamma} 	     & J_{c\theta} \\
		..     & ..     & ..     & J_{\gamma \gamma} & J_{\gamma \theta}\\
		..     & ..     & ..     & ..                & J_{\theta \theta}\\
  \end{bmatrix}
\]
whose elements are reported in the \textit{Appendix}. In order to build the confidence interval and hypothesis tests, we use the fact that the asymptotic distribution of $\boldsymbol{\hat{\eta}}$ can be approximated by the multivariate normal distribution, $N_{6}(\textbf{0},\left[\textbf{J}(\boldsymbol{\hat{\eta}})\right]^{-1} )$, where $\left[\textbf{J}(\boldsymbol{\hat{\eta}})\right]^{-1}$ is the inverse of observed information matrix evaluated at $\boldsymbol{\hat{\eta}}$.

	\section{Application}
In this subsection we consider a real data set concerning the survival times (in days) of guinea pigs injected with different doses of tubercle bacilli. These data, divided by 1000, have already been used in \citet{Kundu2010} in order to illustrate some results derived for the \textit{IWei} distribution.

We compare the fit of the \textit{IWei} model, resulting from the maximum likelihood estimates of the parameters reported by the authors for the uncensored data case, with that obtained by \textit{rGB1-IWei} model.

For the \textit{rGB1-IWe} model, we obtain the following ML estimates and the corresponding standard errors: $\hat{\gamma}=0.8176 \ \ (5.421)$, $\hat{\theta}=0.1284\ \ (0.627)$, $\hat{a}=21.0134 \ \ (192.611)$,  $\hat{b}=76.0581\ \ (720.825)$, $\hat{c}=3.9858\ \ (50.724)$. 
\begin{figure}[ht]
	\centering
		\includegraphics[scale=.6] {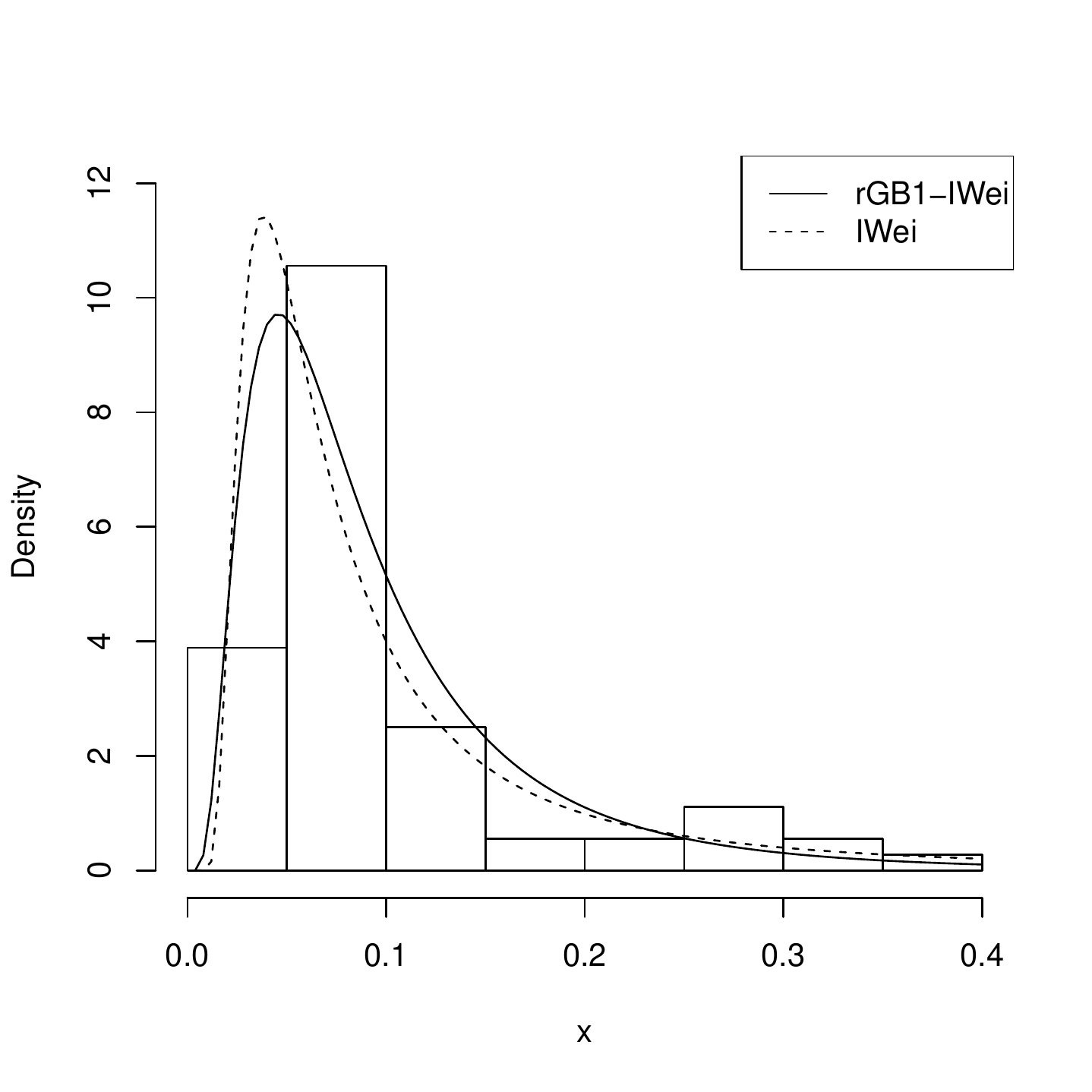}
	\caption{Empirical and fitted densities  of the \textit{rGB1-IWei} and \textit{IWei} distributions.}
	\label{fig.:dens_pigs}
\end{figure}
\begin{figure}[ht]
	\centering
		\includegraphics[scale=.6] {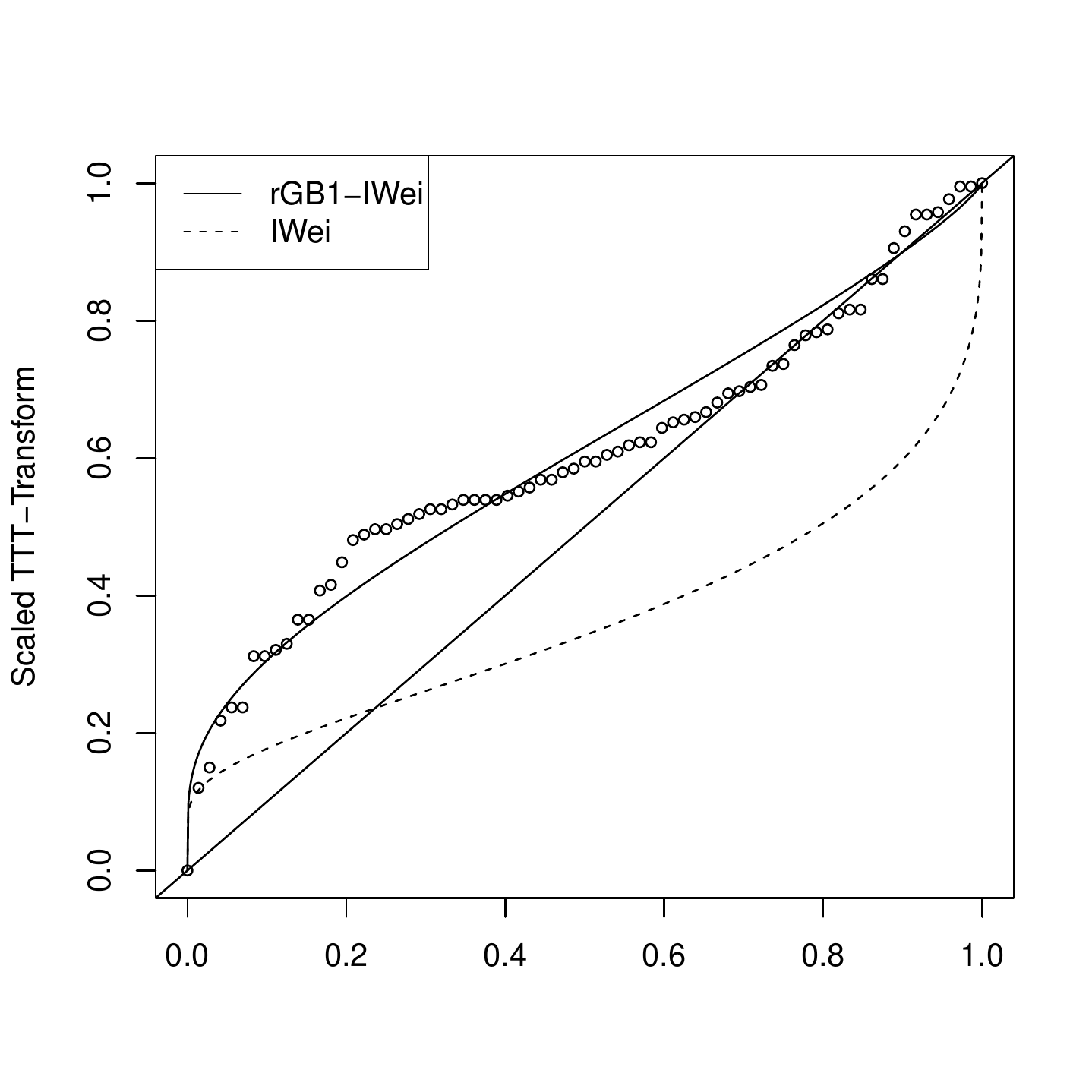}
	\caption{Empirical and fitted Scaled TTT transformation of the Guinea pigs data.}
	\label{fig.:TTT}
\end{figure}
Both the Akaike Information Criterion (AIC) and the  Kolmogorov-Smirnov (KS) statistics suggest that the \textit{rGB1-IWei} model provides a better representation of the data than the \textit{IWei} model. Indeed, considering the different number of parameters and the log-likelihood values reported for the \textit{IWei} and \textit{rGB1-IWei} models ($\hat{l}_{IWei}=101.644$, $\hat{l}_{rGB1-IWei}=107.1106$), we obtain the AIC values, equal to $-199.2888$ and $-204.2213$, respectively. The values of the KS statistic for the two models are $KS_{IWei}=0.1364$ and $KS_{rGB1-IWei}=0.0999$. The associated p-value for the \textit{rGB1-IWei}, equal to 0.4692, supports the hypothesis that the data are drawn from the proposed distribution.
Furthermore, since the \textit{IWei} is a sub-model of the \textit{rGB1-IWei} when $a=b=c=1$, we can consider the maximum values of the unrestricted and restricted log-likelihoods in order to obtain the LR statistics for testing the need of the extra parameters. Based on the LR statistic, equal to 10.932, and the corresponding p-value, equal to 0.0121, at 5\% significance level, we reject the null hypothesis in favor of the new distribution.  The plot of the empirical and the fitted densities are given in Fig. \ref{fig.:dens_pigs}.
Finally, we examine the empirical and the fitted hazard rate. In their paper, \citet{Kundu2010} affirm that the use of the \textit{IWei} model seems reasonable since the plot of the empirical version of the scaled TTT transform for the data indicates that the hazard rate is unimodal. As suggested by \citet{Bergman1984}, we use this plotting procedure as a complementary technique for model identification. If we compare the fitted curves for the \textit{IWei} and the \textit{rGB1-IWei} models (see Fig. \ref{fig.:TTT}), the latter better illustrates the empirical behaviour of the hazard rate, again confirming the superiority of the proposed model in describing such data.

\section{Final Remarks}
In this paper, we introduced the Reflected Generalized Beta of Inverse Weibull Distribution defined by using the Generalized Beta of first type distribution as a generator function, as proposed by \citet{Alexander2012}, considering the reflected version for this generator. This method allows us to obtain a much more flexible distribution than the Inverse Weibull one. Indeed, both the density and the hazard function show additional behaviours compared to the Inverse Weibull distribution. This considerable flexibility makes the new distribution a suitable model for different applications, such as, for example, reliability and survival analysis. 
We provided full treatment of mathematical properties, obtaining moments, entropy, order statistics and a reliability measure. We discussed maximum likelihood estimation and obtained the observed Fisher Information Matrix. Finally, an application to a real data set is given to illustrate the usefulness of the proposed distribution.

\newpage
\section{Appendix}
The elements of the observed information matrix $\boldsymbol {J}(\boldsymbol{\eta})$ for the parameters $(a,b,c,\gamma,\theta)$ are given by

\begin{eqnarray*}
J_{aa}=-n \frac{\ddot{B}_{aa}(a,b)B(a,b)-\left[\dot{B}_{a}(a,b)\right]^2}{\left[B(a,b)\right]^2}
\hspace{0.5 cm}; \hspace{1 cm}
J_{ab}=-n \frac{\ddot{B}_{ab}(a,b)B(a,b)-\dot{B}_{a}(a,b) \dot{B}_{b}(a,b) }{\left[B(a,b)\right]^2}
\end{eqnarray*}

\begin{eqnarray*}
J_{ac}= \sum ^{n}_{i=1} \log \left[ 1- G(x_{i}; \boldsymbol{\tau})\right]
\hspace{0.2 cm}; \hspace{0.3 cm}
J_{bb}=-n \frac{\ddot{B}_{bb}(a,b)B(a,b)-\left[\dot{B}_{b}(a,b)\right]^2}{\left[B(a,b)\right]^2}
\end{eqnarray*}

\begin{eqnarray*}
J_{bc}=-\sum ^{n}_{i=1} \frac{ \left[ 1- G(x_{i}; \boldsymbol{\tau})\right]^c \log \left[ 1- G(x_{i}; \boldsymbol{\tau})\right]  }      {\left( 1-\left[ 1- G(x_{i}; \boldsymbol{\tau})\right]^c \right)}
\end{eqnarray*}

\begin{eqnarray*}
J_{cc}=-\frac{n}{c^2}-(b-1) \sum ^{n}_{i=1} \frac{ \left[ 1- G(x_{i}; \boldsymbol{\tau})\right]^c
\left(\log \left[ 1- G(x_{i}; \boldsymbol{\tau})\right]\right)^2 }
{\left( 1-\left[ 1- G(x_{i}; \boldsymbol{\tau})\right]^c \right)^2}
\end{eqnarray*}

\begin{eqnarray*}
J_{a \tau_{j}}= -c \sum ^{n}_{i=1} \frac{ \dot{G}_{\tau_{j}}(x_{i}; \boldsymbol{\tau})  }{\left[ 1- G(x_{i}; \boldsymbol{\tau})\right]}
\hspace{0.5 cm}; \hspace{0.5 cm}
J_{b\tau_{j}}= c \sum ^{n}_{i=1} \frac{\left[ 1- G(x_{i}; \boldsymbol{\tau})\right]^{c-1} \dot{G}_{\tau_{j}}(x_{i}; \boldsymbol{\tau})   }{\left( 1-\left[ 1- G(x_{i}; \boldsymbol{\tau})\right]^c \right)}
\end{eqnarray*}

\begin{eqnarray*}
J_{c\tau_{j}}=-\sum ^{n}_{i=1} \frac{ \dot{G}_{\tau_{j}}(x_{i}; \boldsymbol{\tau})    }    { \left[ 1- G(x_{i}; \boldsymbol{\tau})\right]  }  \left\{ 1 + (b-1) \frac{  \left[ 1- G(x_{i}; \boldsymbol{\tau})\right]^c
\left( \left[ 1- G(x_{i}; \boldsymbol{\tau})\right]^c -c \log \left[ 1- G(x_{i}; \boldsymbol{\tau})\right] -1 \right) }
{\left( 1-\left[ 1- G(x_{i}; \boldsymbol{\tau})\right]^c \right)^2}
\right\}
\end{eqnarray*}

for $j=1,2$ and $ \tau_1=\gamma$ and $\tau_2=\theta$;

\begin{eqnarray*}
J_{\tau_{j} \tau_{h}}
&=&(1-ac) \sum ^{n}_{i=1} \left( \frac{\ddot{G}_{\tau_{j} \tau_{h}}(x_{i}; \boldsymbol{\tau})\left[ 1- G(x_{i}; \boldsymbol{\tau})\right] + \dot{G}_{\tau_{j}}(x_{i}; \boldsymbol{\tau}) \dot{G}_{\tau_{h}}(x_{i}; \boldsymbol{\tau})  }
{\left[ 1- G(x_{i}; \boldsymbol{\tau})\right]^2 } + \right.\\
&&\left.  \frac{\ddot{g}_{\tau_{j} \tau_{h}}(x_{i}; \boldsymbol{\tau}) g(x_{i}; \boldsymbol{\tau}) -
\dot{g}_{\tau_{j}}(x_{i}; \boldsymbol{\tau}) \dot{g}_{\tau_{h}}(x_{i}; \boldsymbol{\tau})  }
{ \left[g(x_{i}; \boldsymbol{\tau})\right]^2 }\right)
+c(b-1)\sum ^{n}_{i=1} \frac{\left[ 1- G(x_{i}; \boldsymbol{\tau})\right]^{c-2}}{\left( 1-\left[ 1- G(x_{i}; \boldsymbol{\tau})\right]^c \right)^2} \\
&& \left\{
\left( 1-\left[ 1- G(x_{i}; \boldsymbol{\tau})\right]^c \right)
  \left( \dot{G}_{\tau_{j}}(x_{i}; \boldsymbol{\tau}) \dot{G}_{\tau_{h}}(x_{i}; \boldsymbol{\tau})	
	+\left[ 1- G(x_{i}; \boldsymbol{\tau})\right]
\ddot{G}_{\tau_{j}\tau_{h}}(x_{i}; \boldsymbol{\tau}) \right) \right. \\
&& \left. -c  \dot{G}_{\tau_{j}}(x_{i}; \boldsymbol{\tau})  \dot{G}_{\tau_{h}}(x_{i}; \boldsymbol{\tau})
\right\}
\end{eqnarray*}

for $j=1,2$ and $h=1,2$ with $j \leq h$, and $ \tau_1=\gamma$ and $\tau_2=\theta$. \\
Where

\begin{eqnarray*}
\dot{G}_{\gamma} (x_{i}; \boldsymbol{\tau})=-x_{i}^{-\theta} G(x_{i}; \boldsymbol{\tau})
\hspace{0.1 cm}; \hspace{0.2 cm}
\dot{G}_{\theta} (x_{i}; \boldsymbol{\tau})= \gamma x_{i}^{-\theta} \log(\gamma x_{i}) G(x_{i}; \boldsymbol{\tau})
\hspace{0.1 cm}; \hspace{0.2 cm}
\ddot{G}_{\gamma \gamma} (x_{i}; \boldsymbol{\tau})=-x_{i}^{-\theta} \dot{G} _{\gamma }(x_{i}; \boldsymbol{\tau})
\end{eqnarray*}

\begin{eqnarray*}
\ddot{G}_{\gamma \theta} (x_{i}; \boldsymbol{\tau})= x_{i}^{-\theta}
\left\{   \log (x_{i}) G(x_{i}; \boldsymbol{\tau}) - \dot{G} _{\theta}(x_{i}; \boldsymbol{\tau})   \right\}
\end{eqnarray*}

\begin{eqnarray*}
\ddot{G}_{\theta \theta} (x_{i}; \boldsymbol{\tau})= \gamma  x_{i}^{-\theta} \log (x_{i})
\left\{  \dot{G} _{\theta}(x_{i}; \boldsymbol{\tau}) - \log (x_{i}) G(x_{i}; \boldsymbol{\tau}) \right\}
\end{eqnarray*}

\begin{eqnarray*}
\dot{g}_{\gamma} (x_{i}; \boldsymbol{\tau})= \theta x_{i}^{-\theta-1} \left\{ G(x_{i}; \boldsymbol{\tau})+\gamma  \dot{G} _{\gamma}(x_{i}; \boldsymbol{\tau})\right\}
\end{eqnarray*}

\begin{eqnarray*}
\dot{g}_{\theta} (x_{i}; \boldsymbol{\tau})= \gamma \theta x_{i}^{-\theta-1} \left\{ \left[1-\theta \log (x_{i}) \right] G(x_{i}; \boldsymbol{\tau}) + \theta \dot{G} _{\theta}(x_{i}; \boldsymbol{\tau}) \right\}
\end{eqnarray*}

\begin{eqnarray*}
\ddot{g}_{\gamma \gamma} (x_{i}; \boldsymbol{\tau})= \theta x_{i}^{-\theta-1} \left\{
2 \dot{G} _{\gamma}(x_{i}; \boldsymbol{\tau})+\gamma  \ddot{G} _{\gamma \gamma}(x_{i}; \boldsymbol{\tau}) \right\}
\end{eqnarray*}

\begin{eqnarray*}
\ddot{g}_{\gamma \theta} (x_{i}; \boldsymbol{\tau})=  x_{i}^{-\theta-1} \left\{
\left[1-\theta \log (x_{i}) \right] \left[ G(x_{i}; \boldsymbol{\tau})+\gamma \dot{G} _{\gamma}(x_{i}; \boldsymbol{\tau}) \right] + \theta \left[ \dot{G} _{\theta}(x_{i}; \boldsymbol{\tau})
+\gamma  \ddot{G} _{\gamma \theta}(x_{i}; \boldsymbol{\tau}) \right]
\right\}
\end{eqnarray*}

\begin{eqnarray*}
\ddot{g}_{\theta \theta} (x_{i}; \boldsymbol{\tau})= \gamma  x_{i}^{-\theta-1} \left\{
\theta \ddot{G} _{\theta \theta}(x_{i}; \boldsymbol{\tau}) -
\left[ 2-\theta  \log (x_{i})   \right] \log (x_{i}) G(x_{i}; \boldsymbol{\tau}) + \right. \\
\left. 2 \left[ 1-\theta  \log (x_{i})   \right] \dot{G} _{\theta}(x_{i}; \boldsymbol{\tau})
 \right\}
\end{eqnarray*}

\pagebreak

\baselineskip 0.5cm

\addcontentsline{toc}{chapter}{References}
\bibliography{biblio}
\bibliographystyle{plainnat}

\end{document}